# Canonical desingularization in characteristic zero by blowing up the maximum strata of a local invariant

**Edward Bierstone and Pierre D. Milman**[*]

Department of Mathematics, University of Toronto, Toronto, Ontario, Canada M5S 1A1
e-mail: bierston@math.toronto.edu ; milman@math.toronto.edu

Dedicated to Shreeram S. Abhyankar and Heisuke Hironaka

## Contents

[*] The authors' research was partially supported by NSERC grants OGP0009070 and OGP0008949.

## Chapter I. Introduction

This article contains an elementary constructive proof of resolution of singularities in characteristic zero. $\underline{k}$ will denote a field of characteristic zero throughout the paper. Our proof applies to a scheme $X$ of finite type over $\underline{k}$, or to an analytic space $X$ over $\underline{k}$ (in the case that $\underline{k}$ has a complete valuation); we recover, in particular, the great theorems of Hironaka [H1,2], [AHV1,2]. But our work neither was conceived nor is written in the modern language of algebraic geometry. We introduce a discrete local invariant $\text{inv}_X(a)$ whose maximum locus determines a centre of blowing up, leading to desingularization. To define $\text{inv}_X$, we need only to work with a category $\mathcal{A}$ of local-ringed spaces $X = (|X|, \mathcal{O}_X)$ over $\underline{k}$ satisfying the conditions (0.1) following (although further restrictions on $\mathcal{A}$ will be needed for global resolution of singularities).

(0.1) We require that, for each $a \in |X|$: (1) The natural homomorphism $\mathcal{O}_{X,a} \to \widehat{\mathcal{O}}_{X,a}$ into the completion $\widehat{\mathcal{O}}_{X,a} := \varprojlim \mathcal{O}_{X,a}/\underline{m}_{X,a}^{k+1}$ is injective. ($\underline{m}_{X,a}$ denotes the maximal ideal of $\mathcal{O}_{X,a}$.) (2) The residue field $\mathbb{F}_a := \mathcal{O}_{X,a}/\underline{m}_{X,a}$ is included in $\widehat{\mathcal{O}}_{X,a}$. The latter condition implies that there is a surjective homomorphism $\mathbb{F}_a[[X_1, \ldots, X_n]] \to \widehat{\mathcal{O}}_{X,a}$, where $n = \dim_{\mathbb{F}_a} \underline{m}_{X,a}/\underline{m}_{X,a}^2$ and $\mathbb{F}_a[[X_1, \ldots, X_n]]$ denotes the ring of formal power series in the indeterminates $X_1, \ldots, X_n$ with coefficients in $\mathbb{F}_a$.

If $a \in |X|$, then $\text{inv}_X(a)$ depends only on $\widehat{\mathcal{O}}_{X,a}$. More generally, $\text{inv}_X(\cdot)$ is defined recursively over a sequence of formal local blowings-up, $\widehat{\mathcal{O}}_{X,a} = \widehat{\mathcal{O}}_{X_0, a_0} \xrightarrow{\sigma_1^*} \widehat{\mathcal{O}}_{X_1, a_1} \to \cdots \xrightarrow{\sigma_j^*} \widehat{\mathcal{O}}_{X_j, a_j} \to \cdots$, which are "admissible" (in the sense of (1.2) below). For each $j$, $\text{inv}_X(a_j)$ depends only on $\widehat{\mathcal{O}}_{X_j, a_j}$ and on certain information about the previous history of blowings-up which is recorded by the accumulating exceptional divisors. There is an ideal $\widehat{\mathcal{I}}_{S_j}$ of $\widehat{\mathcal{O}}_{X_j, a_j}$ corresponding to a formal "infinitesimal locus" $S_j = S_{\text{inv}_X}(a_j)$ of points $x \in |X_j|$ such that $\text{inv}_X(x) = \text{inv}_X(a_j)$; $S_j$ has only normal crossings. If we choose any component of $S_j$ as the centre of $\sigma_{j+1}^*$, successively for $j = 0, 1, \ldots$, then we get an admissible sequence of formal local blowings-up leading to desingularization of $\widehat{\mathcal{O}}_{X,a}$. In order that the algorithm apply to $\mathcal{O}_{X,a}$, we need to impose conditions on our category $\mathcal{A}$ to guarantee that $\widehat{\mathcal{I}}_{S_j}$ is generated by an ideal $\mathcal{I}_{S_j} \subset \mathcal{O}_{X_j, a_j}$, and $S_j = V(\mathcal{I}_{S_j})$, where $V(\mathcal{I}_{S_j}) = \{x \in |X_j| : f(x) = 0, \text{ for all } f \in \mathcal{I}_{S_j}\}$ (as a germ at $a$).

To obtain a global desingularization algorithm, we need to further restrict $\mathcal{A}$ so that $\text{inv}_X$, defined over an admissible sequence of blowings-up $\cdots \to X_j \xrightarrow{\sigma_j} \cdots \to X_1 \xrightarrow{\sigma_1} X_0$ of $X = X_0$, takes only finitely many maximal values on each $X_j$ (at least locally), and its maximum locus coincides in germs at any point $a_j$ with $S_{\text{inv}_X}(a_j)$ as above. $\text{inv}_X$ then has the property that each local component of its maximum locus extends to a global smooth subspace, justifying the philosophy that "a sufficiently good local choice [of centre of blowing-up] should globalize automatically" [BM4, p. 801].

(0.2) Our desingularization algorithm applies to the following general classes of spaces:

(1) *Algebraic.* Schemes of finite type over $\underline{k}$ (cf. [H1]). Algebraic spaces over $\underline{k}$ (in the sense of Artin [Ar], Knutson [Kn]). Restrictions of schemes $X$ of finite type over $\underline{k}$ to their $\underline{k}$-rational points $|X|_{\underline{k}}$. (Such spaces might be the natural object of study when our main interest lies in the $\underline{k}$-rational points; for instance in resolving the singularities of a real algebraic variety.)

(2) *Analytic.* Real or complex analytic spaces (cf. [H2], [AHV1,2]) or, more generally, analytic spaces in the sense of Serre [Se]. $p$-adic analytic spaces in the sense of Berkovich [Be].

As a rather different example (intermediate between analytic and $\mathcal{C}^\infty$), we can also consider:

(3) "Quasianalytic hypersurfaces", defined by sheaves of principal ideals, each locally generated by a single quasianalytic function, on quasianalytic manifolds in the sense of E.M. Dyn'kin [D].

In each of the classes of (0.2), a space $X$ is locally a subspace of a manifold, or smooth space, $M = (|M|, \mathcal{O}_M)$. For the purpose of global desingularization, a key property of our category of spaces $\mathcal{A}$ is the following:

(0.3) A manifold $M$ in $\mathcal{A}$ can be covered by "regular coordinate charts" $U$: the coordinates $(x_1, \ldots, x_n)$ on $U$ are "regular functions" on $U$ (i.e., each $x_i \in \mathcal{O}_M(U)$) and the partial derivatives $\partial^{|\alpha|}/\partial x^\alpha = \partial^{\alpha_1 + \cdots + \alpha_n}/\partial x_1^{\alpha_1} \cdots \partial x_n^{\alpha_n}$ make sense as transformations $\mathcal{O}_M(U) \to \mathcal{O}_M(U)$. Moreover, for each $a \in U$, there is an injective "Taylor series homomorphism" $T_a \colon \mathcal{O}_{M,a} \to \mathbb{F}_a[[X]] = \mathbb{F}_a[[X_1, \ldots, X_n]]$ such that $T_a$ induces an isomorphism $\widehat{\mathcal{O}}_{M,a} \xrightarrow{\cong} \mathbb{F}_a[[X]]$ and $T_a$ commutes with differentiation: $T_a \circ (\partial^{|\alpha|}/\partial x^\alpha) = (\partial^{|\alpha|}/\partial X^\alpha) \circ T_a$, for all $\alpha \in \mathbb{N}^n$.

In §3 below, we will give a more precise list of the properties of our category of spaces $\mathcal{A}$ that we use to prove global desingularization. As an application of our theorem, we show that desingularization (in the hypersurface case) implies Łojasiewicz's inequalities (§2). (These inequalities seem to be new for quasianalytic functions in dimension $> 2$.) We plan to write another article on the desingularization of "quasi-Noetherian spaces", generalizing Pfaffian varieties in the sense of Khovanskii [Kh].

Our results here were announced in [BM6], and extend in a natural way techniques introduced in [BM3] and [BM4]. When we began thinking about this subject more than fifteen years ago, we were motived by a simple desire to understand how to resolve singularities. One of our goals in this article is that the reader understand the desingularization theorem, rather than simply "know" it is true. We believe that the invariant $\mathrm{inv}_X$ is of interest as a local measure of singularity, beyond desingularization itself. We do not treat non-zero characteristic in this article, although the constructions in Chapter III can largely be made characteristic-free (cf. Remark 1.19 below). Mark Spivakovsky has recently announced a proof of desingularization of arbitrary excellent schemes [Sp].

## 1. An invariant for desingularization

Our proof of resolution of singularities is a variation on our proof of local desingularization ("uniformization") [BM3], [BM4]. We introduce a discrete local invariant $\mathrm{inv}_X(a)$ which is defined recursively over a sequence of blowings-up (or local blowings-up), and which completely determines a succeeding (smooth) centre. Let $X$ denote a space (as above) which is embedded in a manifold (i.e., a smooth space) $M$. Consider a sequence of transformations

$$(1.1) \quad \begin{array}{ccccccccc} \longrightarrow & M_{j+1} & \xrightarrow{\sigma_{j+1}} & M_j & \longrightarrow & \cdots & \longrightarrow & M_1 & \xrightarrow{\sigma_1} & M_0 = M \\ & X_{j+1} & & X_j & & & & X_1 & & X_0 = X \\ & E_{j+1} & & E_j & & & & E_1 & & E_0 = \emptyset \end{array}$$

where, for each $j$, $\sigma_{j+1}: M_{j+1} \to M_j$ denotes a blowing-up (or local blowing-up) with smooth centre $C_j \subset M_j$, $X_{j+1}$ is the strict transform of $X_j$ by $\sigma_j$ (see §3) and $E_{j+1}$ denotes the set of exceptional hypersurfaces. (By definition, $E_{j+1}$ is the set of strict transforms of all $H \in E_j$, together with $\sigma_{j+1}^{-1}(C_j)$. When convenient, we will also use $E_j$ to denote the union of the hypersurfaces $H$ in $E_j$.) If $a \in M_j$, we set $E(a) = \{H \in E_j : a \in H\}$.

Roughly speaking, the goal of "embedded resolution of singularities" is to find a finite sequence of blowings-up (1.1) (or a locally finite sequence in the case of noncompact analytic spaces) such that: If $X'$ and $E'$ denote the final strict transform of $X$ and the final exceptional set (respectively), and if $\sigma: M' \to M$ denotes the composite of the sequence of blowings-up, then (1) $X'$ is smooth; (2) $E' = \sigma^{-1}(\text{Sing } X)$ ($\sigma$ is an isomorphism outside $E'$); (3) $X'$ and $E'$ simultaneously have only normal crossings.

Sing $X$ means the set of singular points of $X$. The condition (3) means that every point of $M'$ admits a coordinate neighbourhood in which $X'$ is a coordinate subspace and each hypersurface $H \in E'$ is a coordinate hypersurface.

Consider a tower of transformations (1.1). Our invariant $\text{inv}_X(a)$, $a \in M_j$, $j = 0, 1, \ldots$, will be defined by induction on $j$, provided that the centres $C_i$, $i < j$, are *admissible* (or $\text{inv}_X$-*admissible*) in the sense that:

(1.2)  (1) $C_i$ and $E_i$ simultaneously have only normal crossings;
        (2) $\text{inv}_X(\cdot)$ is locally constant on $C_i$.

The condition (1.2) (1) guarantees that $E_{i+1}$ is a collection of smooth hypersurfaces having only normal crossings. The notation $\text{inv}_X(a)$, where $a \in M_j$, indicates a dependence on the original space $X$ and not merely on $X_j$. In fact, $\text{inv}_X(a)$, $a \in M_j$, will be an invariant of the local isomorphism class at $a$ of $X_j$ and certain subcollections of $E(a)$ which encode the history of the resolution process. We can think of the desingularization algorithm in the following way: $X \subset M$ determines $\text{inv}_X(a)$, $a \in M$, and thus the first admissible centre of blowing-up $C = C_0$; then $\text{inv}_X(a)$, $a \in M_1$, is defined and determines an admissible blowing-up $\sigma_2$, etc. The exceptional hypersurfaces serve as global coordinate subspaces.

We can allow certain options in the definition of $\text{inv}_X$, but at this point we fix one definition in order to be concrete. $\text{inv}_X(a)$, $a \in M_j$, will be a "word",

(1.3)  $$\text{inv}_X(a) = \left(H_{X_j,a}, s_1(a); \nu_2(a), s_2(a); \ldots, s_t(a); \nu_{t+1}(a)\right),$$

beginning with the *Hilbert-Samuel function* $H_{X_j,a}$ of $X_j$ at $a$; $H_{X_j,a}: \mathbb{N} \to \mathbb{N}$ is the function

$$H_{X_j,a}(\ell) = \dim_{\underline{k}} \frac{\mathcal{O}_{X_j,a}}{\underline{m}_{X_j,a}^{\ell+1}},$$

where $\underline{m}_{X_j,a}$ denotes the maximal ideal of $\mathcal{O}_{X_j,a}$. (In the case of schemes, we would replace $\dim_{\underline{k}}$ by length or $\dim_{\mathbb{F}_a}$ with respect to any embedding $\mathbb{F}_a \hookrightarrow \mathcal{O}_{X_j,a}/\underline{m}_{X_j,a}^{\ell+1}$, where $\mathbb{F}_a$ denotes the residue field $\mathcal{O}_{X_j,a}/\underline{m}_{X_j,a}$ of $a$.)

*Remarks 1.4.* $H_{X_j,a}(\ell)$, for $\ell$ large enough, coincides with a polynomial in $\ell$ of degree $\dim_a X_j$. (See Corollary 3.20.) $H_{X_j,a}(1) - H_{X_j,a}(0) = e_{X_j,a}$ is the minimal embedding dimension of $X_j$ at $a$. (Thus $a \in \text{Sing } X_j$ if and only if $e_{X_j,a} > \dim_a X_j$.) If $a \notin \text{Sing } X_j$, then $H_{X_j,a}(\ell) = \binom{e+\ell}{e}$ for all $\ell$, where $e = e_{X_j,a} = \dim_a X_j$. If $X_j$ is a hypersurface and $\dim_a M_j = n$, then

$$H_{X_j,a}(\ell) = \begin{cases} \binom{n+\ell}{n}, & \ell < \nu_{X_j,a}, \\ \binom{n+\ell}{n} - \binom{n+\ell-\nu_{X_j,a}}{n}, & \ell \geq \nu_{X_j,a}, \end{cases}$$

where $\nu_{X_j,a}$ denotes the *order* of $X_j$ at $a$. ($\nu_{X_j,a} = \max\{\nu : \mathcal{I}_{X_j,a} \subset \underline{m}_{M_j,a}^\nu\}$.) In the hypersurface case, we can therefore replace $H_{X_j,a}$ in the definition of $\mathrm{inv}_X(a)$ by

$$\nu_1(a) = \nu_{X_j,a}.$$

The entries $s_r(a)$ of $\mathrm{inv}_X(a)$ are nonnegative integers reflecting the history of the accumulating exceptional hypersurfaces, and the $\nu_r(a)$, $r \geq 2$, are "multiplicities" of "higher-order terms" of the defining functions of $X_j$ at $a$; $v_2(a), \ldots, v_t(a)$ are quotients of positive integers whose denominators are bounded in terms of the previous part of $\mathrm{inv}_X(a)$. (More precisely, $e_{r-1}!\nu_r(a) \in \mathbb{N}, r = 2, \ldots, t$, where $e_1$ is the smallest integer $k$ such that $H_{X_j,a}(\ell)$ coincides with a polynomial if $\ell \geq k$, and $e_r = \max\{e_{r-1}!, e_{r-1}!\nu_r(a)\}$.) The final entry $\nu_{t+1}(a) = 0$ or $\infty$, and $t \leq n = \dim_a M_j$. (The successive pairs $(\nu_r(a), s_r(a))$ can be defined inductively using functions of $n - r + 1$ variables, so that $t \leq n$ by exhaustion of variables.) The pairs $(\nu_r, s_r)$ originate in the simple uniformization theorem of [BM3, §4]. ($(d, r)$ is the notation used for $(v_1, s_1)$ in [BM3, §4] and in [BM4, §5]; a construction substantially the same is used recursively here to obtain the sequence of pairs comprising $\mathrm{inv}_X$.)

*Example 1.5.* Let $X \subset \underline{k}^n$ denote the hypersurface $x_1^{d_1} + x_2^{d_2} + \cdots + x_t^{d_t} = 0$, where $2 \leq d_1 \leq d_2 \leq \cdots \leq d_t$, $t \leq n$ (i.e., $X$ is defined by the principal ideal generated by $x_1^{d_1} + \cdots + x_t^{d_t}$). Then

$$\mathrm{inv}_X(0) = \left(d_1, 0; \frac{d_2}{d_1}, 0; \cdots; \frac{d_t}{d_{t-1}}, 0; \infty\right).$$

(This is $\mathrm{inv}_X(0)$ at the origin $0$ of $\underline{k}^n$ in "year zero"; i.e., before any blowings-up.)

The simplest form of our embedded desingularization theorem is the following:

**Theorem 1.6.** *Suppose that $|X|$ is quasi-compact. Then there is a finite sequence of blowings-up (1.1) with smooth $\mathrm{inv}_X$-admissible centers $C_j$ such that:*

*(1) For each $j$, either $C_j \subset \mathrm{Sing}\, X_j$ or $X_j$ is smooth and $C_j \subset X_j \cap E_j$.*

*(2) Let $X'$ and $E'$ denote the final strict transform of $X$ and exceptional set, respectively. Then $X'$ is smooth and $X'$, $E'$ simultaneously have only normal crossings.*

(Following Bourbaki, we use "quasi-compact" to mean that every open covering has a finite subcovering; "compact" means "Hausdorff and quasi-compact".) The conclusion of Theorem 1.6 holds, more generally, for $X|U$, where $U$ is any relatively quasi-compact open subset of $|X|$. If $X$ is a non-compact analytic space (for example, over $\mathbb{R}$ or $\mathbb{C}$; see also Theorem 13.3), Theorem 1.6 holds with a locally finite sequence of blowings-up.

If $\sigma$ denotes the composite of the sequence of blowings-up $\sigma_j$, then of course $E'$ is the critical locus of $\sigma$, and $\sigma(E') = \mathrm{Sing}\, X$.

*Remarks 1.7.* (1) Our proof of Theorem 1.6 requires the hypotheses that, for $X$ in our class of spaces $\mathcal{A}$, $\mathrm{Sing}\, X$ is closed and $H_{X,\cdot}$ is upper-semicontinuous, both with respect to the Zariski topology of $|X|$. (The *Zariski topology* of $|X|$ is the topology whose closed sets are of the form $|Y|$, for any closed subspace $Y$ of $X$; see 3.9.) Semicontinuity of $H_{X,\cdot}$ is established in Chapter III (Theorem 9.2), for $X$ in any of the classes of (0.2) (1), (2); in these classes, $\mathcal{O}_X$ is a coherent sheaf of rings, and it follows that $\text{Sing } X$ is Zariski-closed (Proposition 10.1). Both hypotheses above are clear in the hypersurface case, for all classes of (0.2).

(2) Theorem 1.6 resolves the singularities of $X$ in a meaningful geometric sense provided that $\text{Reg } X := |X| \backslash \text{Sing } X$ is Zariski-dense in $|X|$. We will say that $X$ is a *geometric space* if $\text{Reg } X$ is Zariski-dense in $|X|$. For example, *reduced* complex analytic spaces or schemes of finite type are geometric. Suppose that $X$ is a geometric space. If $\sigma: M' \to M$ is a blowing-up with smooth centre $C \hookrightarrow M$, we define the *geometric strict transform* $X''$ of $X$ by $\sigma$ as the smallest closed subspace $Z$ of $\sigma^{-1}(X)$ such that $|Z| \supset |\sigma^{-1}(X)| \backslash |\sigma^{-1}(C)|$. (The strict and geometric strict transforms coincide for reduced schemes or complex analytic spaces.) We can reformulate Theorem 1.6 to resolve the singularities of $X$ by transformations which preserve the class of geometric spaces (§10): we can replace "strict transform" by "geometric strict transform" throughout the desingularization algorithm because, if $X'$ denotes the strict transform of $X$, then $X'' \subset X'$ and, if $a' \in X''$, then $H_{X'',a'} \leq H_{X',a'}$ with equality if and only if $X'' = X'$ in germs at $a'$ (cf. Theorem 1.14 below).

(3) In the categories of (0.2) (1) and (2), algebraic techniques make it possible to use our desingularization algorithm to prove theorems more precise than 1.6; for example, for spaces that are not necessarily reduced (§11). Theorem 1.6 does not exclude the possibility of blowing-up "resolved points"; i.e., a centre of blowing-up $C_j$ prescribed by the desingularization algorithm may include points where $X_j$ is smooth and has only normal crossings with respect to $E_j$. (See Example 2.3.) It is possible to modify $\text{inv}_X$ to avoid blowing up resolved points; see §12.

Our desingularization theorems are presented in Chapter IV. (To be brief, we concentrate in this introduction on an embedded space $X \hookrightarrow M$; see Theorem 13.2 for universal "embedded resolution of singularities" of an abstract space $X$.) We will give a constructive definition of $\text{inv}_X$ in Chapter II. (The main idea will be presented later in this introduction.)

*Remark 1.8. Transforming an ideal to normal crossings* (cf. [H1, Main Th. II]); another example of the principle behind our definition of $\text{inv}_X$: Suppose that $\mathcal{I} \subset \mathcal{O}_M$ is a sheaf of ideals of finite type. Let $\nu_1(a)$ denote the *order* $\nu_{\mathcal{I},a}$ of $\mathcal{I}$ at a point $a \in M$. ($\nu_{\mathcal{I},a} := \max\{\nu : \mathcal{I}_a \subset \underline{m}_{M,a}^\nu\}$.) If $\sigma : M' \to M$ is a local blowing-up with smooth centre $C$, we can define a *weak transform* $\mathcal{I}' \subset \mathcal{O}_{M'}$ of $\mathcal{I}$ by $\sigma$ as follows: For all $a' \in M'$, $\mathcal{I}'_{a'}$ is the ideal generated by $y_{\text{exc}}^{-\nu} f \circ \sigma$, $f \in \mathcal{I}_{\sigma(a')}$, where $\nu$ denotes the generic value of $\nu_1(a)$ on $C$ (and $y_{\text{exc}}$ is a local generator of the ideal of $\sigma^{-1}(C)$ at $a'$). In this context, our construction can be used to extend $\text{inv}_{1/2}(\cdot) = \nu_1(\cdot)$ to an invariant $\text{inv}_{\mathcal{I}}(\cdot)$ which is defined inductively over a sequence of transformations

$$(1.9) \quad \begin{array}{ccccccccc} \longrightarrow & M_{j+1} & \xrightarrow{\sigma_{j+1}} & M_j & \longrightarrow & \cdots & \longrightarrow & M_1 & \xrightarrow{\sigma_1} & M_0 = M \\ & E_{j+1} & & E_j & & & & E_1 & & E_0 = \emptyset \\ & \mathcal{I}_{j+1} & & \mathcal{I}_j & & & & \mathcal{I}_1 & & \mathcal{I}_0 = \mathcal{I} \end{array}$$

where the $\sigma_{j+1}$ are local blowings-up whose successive centres are $\text{inv}_{\mathcal{I}}$-*admissible* (cf. (1.2)), $E_{j+1}$ is the set of exceptional hypersurfaces, and each $\mathcal{I}_{j+1}$ denotes the weak transform of $\mathcal{I}_j$. (See Remark 1.18 below.) Using $\text{inv}_{\mathcal{I}}$, our algorithm gives the following theorem (which is a consequence of Theorem 1.6 in the case that $\mathcal{I} = \mathcal{I}_X$ is the ideal sheaf of a hypersurface $X$).

**Theorem 1.10.** *Suppose that $|M|$ is quasi-compact. Then there is a finite sequence (1.9) of blowings-up $\sigma_j$, $j = 1, \ldots, k$, with smooth $\mathrm{inv}_\mathcal{I}$-admissible centres, such that $\mathcal{I}_k = \mathcal{O}_{M_k}$ and $\sigma^{-1}(\mathcal{I}) := \sigma^*(\mathcal{I}) \cdot \mathcal{O}_{M_k}$ is a normal-crossings divisor, where $\sigma : M_k \to M$ denotes the composite of the $\sigma_j$. ("Normal-crossings divisor" means a principal ideal of finite type, generated locally by a monomial in suitable coordinates.)*

It follows that if $\mathcal{J}_\sigma \subset \mathcal{O}_{M_k}$ denotes the ideal generated (locally, with respect to any coordinate system) by the Jacobian determinant of $\sigma$, then $\mathcal{J}_\sigma \cdot \sigma^{-1}(\mathcal{I})$ is a normal-crossings divisor.

*Remark 1.11.* *In year zero,* there is a straightforward geometric definition of $\mathrm{inv}_X$: Assume that $X$ is a hypersurface. (The following construction will be extended to the general case in Remark 3.23, using the "diagram of initial exponents".) Locally, $X$ is defined by a single equation $f = 0$. Consider the Taylor expansion $f(x) = \sum_{\alpha \in \mathbb{N}^n} f_\alpha x^\alpha$ of $f$ at a point $a$, *for a given coordinate system* $x = (x_1, \ldots, x_n)$. (Say $x(a) = 0$. If $\alpha \in \mathbb{N}^n$, then $x^\alpha$ denotes the monomial $x_1^{\alpha_1} \cdots x_n^{\alpha_n}$. $\mathbb{N}$ always denotes the non-negative integers.) We associate to the Taylor expansion of $f$ at $a$ its *Newton diagram* $\mathfrak{Q}(a) = \{\alpha \in \mathbb{N}^n : f_\alpha \neq 0\}$. Now, let us order the hyperplanes $H$ in $\mathbb{R}^n$ lexicographically with respect to $d = (d_1, \ldots, d_n)$, where the $d_i$ are the intersections of $H$ with the coordinate axes, listed so that $d_1 \leq d_2 \leq \cdots \leq d_n \leq \infty$. We regard $\mathfrak{Q}(a)$ as a subset of the positive orthant of $\mathbb{R}^n$, and let $d(x) = (d_1, \ldots, d_n)$ denote the maximum order of a hyperplane $H$ which lies *under* $\mathfrak{Q}(a)$ (in the sense that for each $(\alpha_1, \ldots, \alpha_n) \in \mathfrak{Q}(a)$, there exists $(\beta_1, \ldots, \beta_n) \in H$ such that $\beta_i \leq \alpha_i$ for each $i$); in particular, $0 < d_1 < \infty$. Of course, $d(x)$ depends on the coordinate system $x = (x_1, \ldots, x_n)$. Set

$$d = \sup_{\substack{\text{coordinate} \\ \text{systems } x}} d(x),$$

$d = (d_1, \ldots, d_n)$. Then

$$\mathrm{inv}_X(0) = \left(d_1, 0; \frac{d_2}{d_1}, 0; \cdots; \frac{d_t}{d_{t-1}}, 0; \infty\right),$$

where $d_t$ is the last finite $d_i$. It is natural to ask whether $d = \sup d(x)$ is realized by a particular coordinate system $x$. (In Example 1.5 above, the supremum is realized by the given coordinates.) The construction we use to define $\mathrm{inv}_X$ in Chapter II gives a positive answer. Moreover, beginning with any coordinate system, we find an explicit change of variables to obtain coordinates $x = (x_1, \ldots, x_n)$ satisfying a criterion which guarantees that $d(x) = d$; in these coordinates, the centre of the first blowing-up in our resolution algorithm is $x_1 = \cdots = x_t = 0$ (where the coordinates are indexed so that $d_i$ corresponds to $x_i$, for each $i$). Consider also a second coordinate system $y = (y_1, \ldots, y_n)$ in which the supremum $d$ is realized (indexed again so that $d_i$ corresponds to $y_i$, for each $i$). Write $y = \varphi(x)$, $\varphi = (\varphi_1, \ldots, \varphi_n)$, for the coordinate transformation. Then, using $w_i = d_1/d_i$ as weights for the coordinates $x_i$ and $y_i$, $i = 1, \ldots, n$ (cf. Remark 3.23), the weighted initial forms of $f$ with respect to $x$ and $y$ are obtained one from the other by the substitution $y = \varphi_w(x)$, where each $y_i = \varphi_{w,i}(x)$ denotes the weighted homogeneous part of order $w_i$ in the Taylor expansion of $y_i = \varphi_i(x)$. (This remark will not be used in this article; we intend to pursue it elsewhere.)

*Remark 1.12. A combinatorial analogue of resolution of singularities.* Let $\mathcal{M}$ be a finite simplicial complex. We define the (*simplicial*) *blowing-up* $\sigma_\Sigma$ of $\mathcal{M}$ *along* a simplex $\Sigma$ as the smallest simplicial subdivision $\mathcal{M}'$ of $\mathcal{M}$ which includes the barycentre of $\Sigma$. If $V(\mathcal{M})$ denotes the set of vertices (0-simplices) $\{H_1, \ldots, H_d\}$ of $\mathcal{M}$, then $V(\mathcal{M}') = \{H'_1, \ldots, H'_{d+1}\}$, where $H'_k = H_k$, $k \leq d$, and $H'_{d+1}$ is the barycentre of $\Sigma$. If $D$ is a function $D: V(\mathcal{M}) \to \mathbb{Z}$, we define the transform $D'$ of $D$ by $\sigma_\Sigma$ as $D': V(\mathcal{M}') \to \mathbb{Z}$, where $D'(H'_k) = D(H_k)$, $k \leq d$, and $D'(H'_{d+1}) = \sum_{H_k \in V(\Sigma)} D(H_k)$. If $D_1, D_2: V(\mathcal{M}) \to \mathbb{Z}$, we say that $D_1 \leq D_2$ if $D_1(H_k) \leq D_2(H_k)$ for all $k$. *Theorem.* Suppose that $D_j: V(\mathcal{M}) \to \mathbb{Z}$, $j = 1, \ldots, s$. Then there is a finite sequence of simplicial blowings-up of $\mathcal{M}$ after which the transforms $D'_j$ of the $D_j$ are locally totally ordered in the following sense: Let $\mathcal{M}'$ denote the final transform of $\mathcal{M}$. Then, for every simplex $\Sigma$ of $\mathcal{M}'$, there is a permutation $(j_1, \ldots, j_s)$ of the indices $j$ such that $D'_{j_1}|V(\Sigma) \leq D'_{j_2}|V(\Sigma) \leq \cdots \leq D'_{j_s}|V(\Sigma)$.

A finite simplicial complex can be associated to a system of smooth hypersurfaces with only normal crossings in a smooth ambient space. Let $M$ be a manifold and let $E$ denote a finite collection of smooth hypersurfaces $\{H_1, \ldots, H_d\}$ in $M$ having only normal crossings. We associate to $E$ the simplicial complex $\mathcal{M} = \mathcal{M}(E)$ whose vertices correspond to the hypersurfaces $H_k$ and whose simplices $\Sigma$ correspond to nonempty intersections $H_{k_1} \cap \cdots \cap H_{k_q}$. Every finite simplicial complex can be realized in this way. (This follows from the preceding theorem, for example.) Let us say that a blowing-up $\sigma: M' \to M$ is *admissible* if its centre $C$ is an intersection of hypersurfaces in $E$; i.e., $C = C(\Sigma)$ corresponds to a simplex $\Sigma$ of $\mathcal{M}$. The system of hypersurfaces $E$ transforms under $\sigma$ to $E' = \{H'_1, \ldots, H'_{d+1}\}$, where $H'_k$ denotes the strict transform of $H_k$, $k \leq d$, and $H'_{d+1} = \sigma^{-1}(C)$. It is easy to see that $\mathcal{M}' = \mathcal{M}(E')$ is the simplicial blowing-up $\sigma_\Sigma$ of $\mathcal{M}$. A formal divisor $D = \sum_{k=1}^{d} n_k[H_k]$ (where each $n_k \in \mathbb{Z}$) on $M$ corresponds to the function $D(H_k) = n_k$ on $V(\mathcal{M})$. The (*total*) *transform* of $D = \sum n_k[H_k]$ *by* $\sigma$ is defined as $D' = \sum_{k=1}^{d} n_k[H'_k] + \big(\sum_{H_k \in V(\Sigma)} n_k\big)[H'_{d+1}]$. Clearly, this is the same as the combinatorial transformation rule above. The preceding theorem is equivalent to the following "combinatorial desingularization theorem" (which should be compared to the role played by Lemma 4.7 in [BM3, §4]).

**Theorem 1.13.** *Let $M$ be a manifold and let $E = \{H_1, \ldots, H_d\}$ be a finite collection of smooth hypersurfaces in $M$ having only normal crossings. Suppose we have a system of formal divisors $D_j = \sum n_{jk}[H_k]$ (where each $n_{jk} \in \mathbb{Z}$), $j = 1, \ldots, s$. Then there is a finite sequence of admissible blowings-up of $M$ after which the transforms $D'_j$ of the $D_j$ are locally totally ordered in the following sense: Let $M', E'$ denote the final transforms of $M, E$. Then, for each $a' \in M'$, there is a permutation $(j_1, \ldots, j_s)$ of the indices $j$ such that $D'_{j_1}(H) \leq D'_{j_2}(H) \leq \cdots \leq D'_{j_s}(H)$ for all $H \in E'$ such that $a' \in H$.*

*Proof.* It is enough to consider $s = 2$. The following argument is a very simple parallel of the construction in Chapter II (and Theorem 1.14). Let $a \in M$. Set $\nu_1(a) := \min\{\sum_{H \ni a} D_1(H), \sum_{H \ni a} D_2(H)\} - \sum_{H \ni a} \min\{D_1(H), D_2(H)\}$, and $\mathrm{inv}(a) := \nu_1(a)$. ("$s_1(a)$" is not needed and "$\nu_2(a)$" = 0.) Put

$$\mu_2(a) := \max\left\{\sum_{H \ni a} D_1(H), \sum_{H \ni a} D_2(H)\right\} - \sum_{H \ni a} \min\{D_1(H), D_2(H)\} \,.$$

Let $S_a$ denote (the germ at $a$ of) $\{x \in M : \nu_1(x) = \nu_1(a)\}$. Then the irreducible components $Z$ of $S_a$ are of the form $Z = Z_I := \bigcap_{H \in I} H$ for certain $I \subset E$ (as germs at $a$) (and $Z = S_a \cap \bigcap_{H \supset Z} H$; cf. Theorem 1.14 (3)): To be explicit, say that $\nu_1(a) = \sum_{H \ni a}(D_1(H) - D_0(H))$, where $D_0(H) := \min\{D_1(H), D_2(H)\}$. Set $J_k(a) := \{H : H \ni a$ and $D_k(H) > D_0(H)\}$, $k = 1, 2$. Then each $Z = Z_I$, where $I = J_1(a) \cup J$ and $J$ is a subset of $J_2(a)$ that is minimal with respect to the property that $\sum_{H \in J}(D_2(H) - D_0(H)) \geq \nu_1(a)$. (In particular, if $\mu_2(a) = \nu_1(a)$, then $S_a = Z_I$, where $I = J_1(a) \cup J_2(a)$.)

Now let $\nu_1$ denote the maximum value of the invariant $\mathrm{inv}(a)$, $a \in M$, and let $S := \{x \in M : \nu_1(x) = \nu_1\}$. Then the irreducible components of $S$ are the $Z_I$ above, for all $a \in S$. Write $\mu_2(I) := \min_{a \in Z_I} \mu_2(a)$; then $\mu_2(I) = \max\{\sum_{H \in I}(D_1(H) - D_0(H)), \sum_{H \in I}(D_2(H) - D_0(H))\} \geq \nu_1$. Let $\sigma$ denote the blowing-up with centre one of these components $Z_I$. We claim that $(\nu_1(b); \mu_2(b)) < (\nu_1(\sigma(b)); \mu_2(\sigma(b)))$, for all $b \in \sigma^{-1}(Z_I)$ (so the theorem follows by induction). Indeed, by the minimality property above, $I = J_1 \cup J_2$, where $J_k := \{H \in I : D_k(H) > D_0(H)\}$, $k = 1, 2$. Say that $\nu_1 = \sum_{H \in I}(D_1(H) - D_0(H))$. Let $a \in Z_I$. If $\nu_1 < \mu_2(I)$, then $J_1 = J_1(a)$; if $\nu_1 = \mu_2(I)$, then we can assume that the same is true by interchanging $k = 1$ and 2 if necessary. In any case, $0 \leq \sum_{H \in I}(D_2(H) - D_0(H)) - \nu_1 < D_2(H_*) - D_0(H_*)$ for every $H_* \in J_2$. Let $b \in H_I := \sigma^{-1}(Z_I)$ and let $a = \sigma(b)$. Then $D_1(H_I) := \sum_{H \in I} D_1(H) \leq \sum_{H \in I} D_2(H) =: D_2(H_I)$; hence $D_0(H_I) = D_1(H_I)$, and $D_2(H_I) - D_0(H_I) < D_2(H_*) - D_0(H_*)$ for every $H_* \in J_2$. If $H \in E$, let $H'$ denote the strict transform of $H$. Since $D_1(H_I) - D_0(H_I) = 0$, it follows that $\nu_1(b) \leq \nu_1(a)$ and $\nu_1(b) = \nu_1(a)$ if and only if $b \in \bigcap_{H \in J_1} H'$ and $\mu_2(b) = \sum_{H \ni b}(D_2(H) - D_0(H)) \geq \sum_{H \ni b}(D_1(H) - D_0(H)) = \nu_1(b)$ (in particular, $b \notin H'_*$, for some $H_* \in J_2$). But $\mu_2(a) = \sum_{H \ni a}(D_2(H) - D_0(H))$ since $\nu_1 = \nu_1(a)$ and $J_1 = J_1(a)$. Therefore, $\nu_1(b) = \nu_1(a)$ implies that $\mu_2(b) \leq \mu_2(a) - (D_2(H_*) - D_0(H_*)) + (D_2(H_I) - D_0(H_I)) < \mu_2(a)$, as required. (Since $\nu_1(b) \leq \mu_2(b)$, it follows that $\nu_1(b) < \nu_1(a)$ if $\mu_2(a) = \nu_1(a)$.) □

*Fundamental properties of* $\mathrm{inv}_X$

Let $X$ denote a closed subspace of a smooth space $M$, as before. Our desingularization theorems will follow from four key properties satisfied by $\mathrm{inv}_X$, for any admissible sequence of blowings-up; we list these properties in Theorem 1.14 following.

A function $\tau : |M| \to \Sigma$ with values in a partially-ordered set $\Sigma$ will be called *Zariski-semicontinuous* if $\tau$ locally takes only finitely many values and, for all $\sigma \in \Sigma$, $S_\sigma := \{x \in |M| : \tau(x) \geq \sigma\}$ is Zariski-closed. If $|X|$ is Noetherian, then $\tau$ is Zariski-semicontinuous if and only if each $a \in |X|$ admits a Zariski-open neighbourhood in which $\tau(x) \leq \tau(a)$. (See Lemma 3.10 and Definition 3.11.)

The function $\tau(\cdot) = H_{X,\cdot}$ takes values in the set $\mathbb{N}^{\mathbb{N}}$ of functions from $\mathbb{N}$ to itself. $\mathbb{N}^{\mathbb{N}}$ can be partially ordered as follows: If $H_1, H_2 \in \mathbb{N}^{\mathbb{N}}$, then $H_1 < H_2$ if $H_1(\ell) \leq H_2(\ell)$ for all $\ell$, and $H_1(\ell) < H_2(\ell)$ for some $\ell$. We can then use the lexicographic ordering of sequences of the form (1.3) to obtain a partially-ordered set in which $\mathrm{inv}_X(\cdot)$ takes values.

**Theorem 1.14.** *Consider any* $\mathrm{inv}_X$*-admissible sequence of local blowings-up (1.1). The following properties hold.*

*(1) Semicontinuity. (i) For each $j$, every point of $|M_j|$ admits a neighbourhood $U$ such that $\mathrm{inv}_X$ takes only finitely many values in $U$ and, for all $a \in U$, $\{x \in U : \mathrm{inv}_X(x) \leq \mathrm{inv}_X(a)\}$ is Zariski-open in $|M_j|U|$. (ii) $\mathrm{inv}_X$ is "infinitesimally upper-semicontinuous" in the sense that $\mathrm{inv}_X(a) \leq \mathrm{inv}_X\bigl(\sigma_j(a)\bigr)$ for all $a \in M_j$, $j \geq 1$.*

*(2) Stabilization. Given $a_j \in M_j$ such that $a_j = \sigma_{j+1}(a_{j+1})$, $j = 0, 1, 2, \ldots$, there exists $j_0$ such that $\mathrm{inv}_X(a_j) = \mathrm{inv}_X(a_{j+1})$ when $j \geq j_0$. (In fact, any nonincreasing sequence in the value set of $\mathrm{inv}_X$ stabilizes.)*

*(3) Let $a \in M_j$ and let $S_X(a)$ denote the germ at $a$ (with respect to the Zariski topology) of $S_{\mathrm{inv}_X(a)}$ (so that $\mathrm{inv}_X(\cdot) = \mathrm{inv}_X(a)$ on $S_X(a)$). Then $S_X(a)$ and $E(a)$ simultaneously have only normal crossings. If $\mathrm{inv}_X(a) = (\ldots; \infty)$, then $S_X(a)$ is smooth. If $\mathrm{inv}_X(a) = (\ldots; 0)$ and $Z$ denotes an irreducible component of $S_X(a)$, then*

$$Z = S_X(a) \cap \bigcap \{H \in E(a) : Z \subset H\}.$$

*(4) Let $a \in M_j$. If $\mathrm{inv}_X(a) = (\ldots; \infty)$ and $\sigma$ is the local blowing-up of $M_j$ with centre $S_X(a)$, then $\mathrm{inv}_X(a') < \mathrm{inv}_X(a)$ for all $a' \in \sigma^{-1}(a)$. In general, there is an additional invariant $\mu_X(a) \geq 1$ with the following property: Let $Z$ be an irreducible component of $S_X(a)$ and let $\sigma$ denote the local blowing-up with centre $Z$; then $\bigl(\mathrm{inv}_X(a'), \mu_X(a')\bigr) < \bigl(\mathrm{inv}_X(a), \mu_X(a)\bigr)$ for all $a' \in \sigma^{-1}(a)$. ($e_t! \mu_X(a) \in \mathbb{N}$, with $e_t$ as defined following Remarks 1.4.)*

Theorem 1.14 will be proved in Chapter II in the case that $X$ is a hypersurface, and completed in Chapter III in the general case. Condition (1) (i) implies that $\mathrm{inv}_X$ is Zariski-semicontinuous if $|X|$ is quasi-compact or if $X$ is an analytic space over a locally compact field $\underline{k}$ (Remark 6.14). Note that, because of the bounds on the denominators of the terms $\nu_r(a)$ in $\mathrm{inv}_X(a)$, the stabilization property (2) of Theorem 1.14 is an immediate consequence of the corresponding property of the Hilbert-Samuel function. An elementary proof of stability of the Hilbert-Samuel function can be found in [BM4, Th. 5.2.1]. The present article is self-contained except for this result and some elementary properties of the diagram of initial exponents, for which we give references in §3 below.

*Remarks 1.15.* (1) Let $a \in M_j$, for some $j$, and let $U = \{x \in |M_j| : \mathrm{inv}_X(x) \leq \mathrm{inv}_X(a)\}$. Then each irreducible component of $S_X(a)$ extends to a smooth (Zariski-) closed subset of $U$. This is a consequence of property (3) of Theorem 1.14: If $a \in M_j$, we label every component $Z$ of $S_X(a)$ as $Z_I$, where $I = \{H \in E(a) : Z \subset H\}$. Consider any total ordering on the collection of all subsets $I$ of $E_j$. For each $a \in M_j$, put $J(a) = \max\{I : Z_I \text{ is a component of } S_X(a)\}$; set

$$\mathrm{inv}_X^e(a) = \bigl(\mathrm{inv}_X(a); J(a)\bigr).$$

Clearly, $\mathrm{inv}_X^e(\cdot)$ satisfies property (1) (i) of Theorem 1.14 and its locus of maximal values on $U$ is smooth.

Of course, given $a \in M_j$ and any component $Z_I$ of $S_X(a)$, we can choose the ordering above so that $I = J(a) = \max\{J : J \subset E_j\}$; therefore, $Z_I$ extends to a smooth Zariski-closed subset of $U$.

(2) The preceding construction shows that $\text{inv}_X(a)$ can be extended to an invariant $\text{inv}_X^e(a) = (\text{inv}_X(a); J(a))$ which satisfies properties (1)–(4) of Theorem 1.14 and has the additional property that, for all $a$, $S_X^e(a)$ is smooth (where $S_X^e(a)$ denotes the germ of $S_{\text{inv}_X^e(a)}$ at $a$): It suffices to order the subsets of each $E_j$ in the following way: Write $E_j = \{H_1^j, \ldots, H_j^j\}$, where $H_i^j$ denotes the strict transform of $H_i^{j-1}$ by $\sigma_j$, $i = 1, \ldots, j-1$, and $H_j^j = \sigma_j^{-1}(C_{j-1})$ (i.e., each $H_i^j$ is the strict transform of $\sigma_i^{-1}(C_{i-1})$ by the sequence of blowings-up $\sigma_{i+1}, \ldots, \sigma_j$). Associate to each $I \subset E_j$ the sequence $(\delta_1, \ldots, \delta_j)$, where $\delta_i = 0$ if $H_i^j \notin I$ and $\delta_i = 1$ if $H_i^j \in I$, and use the lexicographic ordering of such sequences, for all $j$ and all $I \subset E_j$.

*Universal and canonical desingularization*

We can use the extended invariant $\text{inv}_X^e$ and Theorem 1.14 to obtain a desingularization algorithm with uniquely determined centres of blowing up: When our spaces are quasi-compact (e.g., in the categories of schemes or compact analytic spaces) we get a tower of $\text{inv}_X$-admissible blowings-up (1.1) by successively choosing as each smooth closed centre $C_j$, the locus of (the finitely many) maximal values of $\text{inv}_X^e$ on $\text{Sing}\, X_j$. (If $a \in \text{Sing}\, X_j$, then $S_X(a) \subset \text{Sing}\, X_j$ because the Hilbert-Samuel function already distinguishes between smooth and singular points.) By property (4) of Theorem 1.14, $(\text{inv}_X(a'), \mu_X(a')) < (\text{inv}_X(a), \mu_X(a))$ for all $a \in C_j$ and $a' \in \sigma_{j+1}^{-1}(a)$. The embedded desingularization theorem 1.6 follows. (See §10 below.) The desingularization theorem for analytic spaces $X$ which are not necessarily compact follows from the algorithm applied to relatively compact open subsets of $X$ (§13).

*Remark 1.16.* It is not difficult to define an extended invariant $\text{inv}_X^e$ with the additional property that, for all $a \in M_j$, if $\sigma$ is a local blowing-up with centre $S_X^e(a)$ (the germ at $a$ of $S_{\text{inv}_X^e(a)}$), then $\text{inv}_X^e(a') < \text{inv}_X^e(a)$ for all $a' \in \sigma^{-1}(a)$ (so that $\mu_X(\cdot)$ is not needed in property (4) of Theorem 1.14 formulated for $\text{inv}_X^e$); see Remark 6.17 below.

Our desingularization algorithm is *universal* for Noetherian spaces: To every $X$, we associate a morphism of resolution of singularities $\sigma_X \colon X' \to X$ such that any local isomorphism $X \,|\, U \to Y \,|\, V$ lifts to an isomorphism $X' \,|\, \sigma_X^{-1}(U) \to Y' \,|\, \sigma_Y^{-1}(V)$ (in fact, lifts to isomorphisms throughout the entire towers of blowings-up). ($U, V$ denote Zariski-open subsets of $|X|, |Y|$, respectively.) See §13.

For analytic spaces which are not necessarily compact, the resulting desingularization procedure is *canonical*: Given $X$, there is a morphism of desingularization $\sigma_X \colon X' \to X$ such that any isomorphism $X|U \to X|V$ (over subsets $U, V$ of $|X|$ which are open in the Hausdorff topology) lifts to an isomorphism $X' \,|\, \sigma_X^{-1}(U) \to X' \,|\, \sigma_X^{-1}(V)$. (See §13.)

*Presentation of the invariant*

We outline here the constructive definition of $\text{inv}_X$ that is detailed in Chapter II below. (It might be helpful to read this subsection in parallel with the examples of §2.) The entries $s_1(a), \nu_2(a), s_2(a), \ldots$ of $\text{inv}_X(a)$ (1.3) will themselves be defined recursively. Let us write $\text{inv}_r$ for $\text{inv}_X$ truncated after the entry $s_r$ (with the convention that $\text{inv}_r(a) = \text{inv}_X(a)$ if $r > t$). We also write $\text{inv}_{r+\frac{1}{2}}(a) = (\text{inv}_r; \nu_{r+1})$, so that $\text{inv}_{1/2}(a)$ means $H_{X_j, a}$. For each $r$, the entries $s_r, \nu_{r+1}$ of $\text{inv}_X$ can, in fact, be defined inductively over a tower of (local) blowings-up (1.1) whose centres $C_i$ are merely $(r - \frac{1}{2})$-admissible in the sense that:

(1.17)(1) $C_i$ and $E_i$ simultaneously have only normal crossings;

(2) $\text{inv}_{r-\frac{1}{2}}$ is locally constant on $C_i$.

Once $\text{inv}_{r+\frac{1}{2}}$ is defined, $s_{r+1}$ can be introduced immediately, in an invariant way: Consider a tower of local blowings-up (1.1) with $(r + \frac{1}{2})$-admissible centres. Write $\pi_{ij} = \sigma_{i+1} \circ \cdots \circ \sigma_j, i = 0, \ldots, j-1$, and $\pi_{jj} = id$. Suppose $a \in M_j$. We set $a_i = \pi_{ij}(a)$. First consider $r = 0$. Let $i$ denote the smallest index $k$ such that $\text{inv}_{1/2}(a) = \text{inv}_{1/2}(a_k)$ and set $E^1(a) = \{H \in E(a) : H \text{ is the strict transform of some hyperplane in } E(a_i)\}$. We define $s_1(a) = \#E^1(a)$. In general, suppose that $i$ is the smallest index $k$ such that $\text{inv}_{r+\frac{1}{2}}(a) = \text{inv}_{r+\frac{1}{2}}(a_k)$. Let $E^{r+1}(a) = \{H \in E(a) \setminus \bigcup_{q \leq r} E^q(a) : H \text{ is the strict transform of some element of } E(a_i)\}$. We define $s_{r+1}(a) = \#E^{r+1}(a)$.

We will introduce each $\nu_{r+1}(a)$ and prove Theorem 1.14 above by an explicit construction in local coordinates. The central idea of our approach to desingularization is that the local construction survives certain blowings-up. This idea (which we use, for example, to prove that the numerical characters $\nu_{r+1}(a)$ are invariant) already appears in our earlier work [BM 3,4,6], but the notion of a local "presentation of an invariant", introduced here, captures it in a better way. (We recommend reading [BM3, §4] for the essence of the approach in the simplest possible context.)

Let us consider data of the following type at a (closed) point $a \in M$ (say that $a$ is $\underline{k}$-rational, in the case of schemes):

$N = N(a)$: a germ at $a$ of a regular submanifold of $M$ of codimension $p$;

$\mathcal{H}(a) = \{(h, \mu_h)\}$: a finite collection of pairs $(h, \mu_h)$, where each $h \in \mathcal{O}_{N,a}$ and each $\mu_h \in \mathbb{Q}$ is an "assigned multiplicity" $\mu_h \leq \mu_a(h)$ (where $\mu_a(h)$ denotes the order of $h$ at $a$);

$\mathcal{E}(a)$: a collection of smooth hypersurfaces $H \ni a$ such that $N$ and $\mathcal{E}(a)$ simultaneously have only normal crossings, and $N \not\subset H$, for all $H \in \mathcal{E}(a)$.

We will call $(N(a), \mathcal{H}(a), \mathcal{E}(a))$ an *infinitesimal presentation*, and we define its *equimultiple locus* $S_{\mathcal{H}(a)}$ as $\{x \in N : \mu_x(h) \geq \mu_h, \text{ for all } (h, \mu_h) \in \mathcal{H}(a)\}$. $S_{\mathcal{H}(a)} \subset N$ is well-defined as a germ at $a$. Given an infinitesimal presentation $(N(a), \mathcal{H}(a), \mathcal{E}(a))$, we also define a transform $(N(a'), \mathcal{H}(a'), \mathcal{E}(a'))$ by a morphism of each of 3 types: (i) *admissible blowing-up*, (ii) *projection from the product with a line*, (iii) *exceptional blowing-up*. See (4.3) below for the definitions. An admissible blowing-up, for example, means a local blowing-up $\sigma: M' \to M$ with smooth centre $C \subset S_{\mathcal{H}(a)}$ such that $C$ and $\mathcal{E}(a)$ simultaneously have only normal crossings. In this case, let $N'$ denote the strict transform of $N$ by $\sigma$, and let $a' \in \sigma^{-1}(a)$ such that $a' \in N'$ and $\mu_{a'}(h') \geq \mu_h$, for all $(h, \mu_h) \in \mathcal{H}(a)$, where $h' = y_{\text{exc}}^{-\mu_h} h \circ \sigma$ (provided that such $a'$ exists). ($y_{\text{exc}}$ de-notes a local generator of the ideal of $\sigma^{-1}(C)$.) We set $N(a')$ = the germ of $N'$ at $a'$, $\mathcal{H}(a') = \{(h', \mu_h)\}$, and $\mathcal{E}(a') = \{\sigma^{-1}(C)\} \cup \{H' : H \in \mathcal{E}(a), a' \in H'\}$ (where $H'$ = the strict transform of $H$).

Transformations of types (ii) and (iii) will be needed only to prove the invariance of $\nu_{r+1}(a)$ using certain sequences of test blowings-up. Of course $y_{\text{exc}}^{-\mu_h} h \circ \sigma$ above is defined only up to an invertible factor, but two different choices are equivalent in the sense of the following definition: Given $\mathcal{E}(a)$, we will say that two infinitesimal presentations $(N, \mathcal{F}(a), \mathcal{E}(a))$ and $(P, \mathcal{H}(a), \mathcal{E}(a))$ (perhaps of different codimension) are *equivalent* (*with respect to transformations of* types (i), (ii) *and* (iii)) if:

(1) $S_{\mathcal{F}(a)} = S_{\mathcal{H}(a)}$.

(2) If $\sigma$ is a local blowing-up as in (i) and $a' \in \sigma^{-1}(a)$, then $a' \in N'$ and $\mu_{a'}(y_{\text{exc}}^{-\mu_f} f \circ \sigma) \geq \mu_f$, for all $(f, \mu_f) \in \mathcal{F}(a)$, if and only if $a' \in P'$ and $\mu_{a'}(y_{\text{exc}}^{-\mu_h} h \circ \sigma) \geq \mu_h$, for all $(h, \mu_h) \in \mathcal{H}(a)$.

(3) After a transformation of type (i), (ii) or (iii), $(N', \mathcal{F}(a'), \mathcal{E}(a'))$ is equivalent to $(P', \mathcal{H}(a'), \mathcal{E}(a'))$. (This makes sense recursively.)

For example, assume that $(N(a), \mathcal{H}(a), \mathcal{E}(a))$ is an infinitesimal presentation, where $\mathcal{H}(a) = \{(h, \mu_h)\}$. Then: (1) There is an equivalent infinitesimal presentation with $\mu_h \in \mathbb{N}$, independent of $h$: we can simply replace each $(h, \mu_h)$ by $(h^{e/\mu_h}, e)$, for suitable $e$. (2) Suppose there is $(h, \mu_h) \in \mathcal{H}(a)$ with $\mu_a(h) = \mu_h$ and $h = \Pi h_i^{m_i}$. If we replace $(h, \mu_h)$ in $\mathcal{H}(a)$ by the collection of pairs $(h_i, \mu_{h_i})$, where each $\mu_{h_i} = \mu_a(h_i)$, then we obtain an equivalent infinitesimal presentation.

We will prove that

$$\mu_{\mathcal{H}(a)} := \min_{(h, \mu_h) \in \mathcal{H}(a)} \frac{\mu_a(h)}{\mu_h}$$

is an invariant of the equivalence class of the infinitesimal presentation $(N(a), \mathcal{H}(a), \mathcal{E}(a))$ (in fact, with respect to transformations of types (i), (ii) alone).

The starting point of our construction is a local invariant which admits a presentation; we consider here the Hilbert-Samuel function $H_{X,\cdot}$ of our space $X \subset M$ (but see also Remarks 1.8, 1.18): We first introduce the transform $X'$ of $X$ by a morphism $\sigma$ of type (i), (ii) or (iii): $X'$ is the strict transform of $X$ in the case of (i), and the total transform $\sigma^{-1}(X)$ in the case of (ii) or (iii). An infinitesimal presentation $N = (N(a), \mathcal{H}(a), \mathcal{E}(a))$ with codim $N = p$ will be called a (*codimension p*) *presentation of* $H_{X,\cdot}$ *at $a$* (*with respect to* $\mathcal{E}(a)$) if:

(1) $S_{\mathcal{H}(a)} = S_H(a)$, where $S_H(a)$ denotes the germ at $a$ of the *Hilbert-Samuel stratum* $\{x : H_{X,x} = H_{X,a}\}$.

(2) After an admissible local blowing-up $\sigma$ ((i) above), $H_{X',a'} = H_{X,a}$ if and only if $a' \in N'$ and $\mu_{a'}(h') \geq \mu_h$ for all $(h, \mu_h) \in \mathcal{H}(a)$).

(3) Conditions (1) and (2) continue to hold after any sequence of transformations of types (i), (ii) and (iii).

In particular, after any sequence of transformations (i), (ii), (iii), the transform $(N(a'), \mathcal{H}(a'), \mathcal{E}(a'))$ is a (codimension $p$) presentation of $H_{X',\cdot}$ at $a'$, with respect to $\mathcal{E}(a')$. Of course, any two presentations of $H_{X,\cdot}$ at $a$ with respect to $\mathcal{E}(a)$ are equivalent. It is clear that the equivalence class of a presentation of $H_{X,\cdot}$ at $a$ with respect to $\mathcal{E}(a)$ depends only on the local isomorphism class of $M$, $X$, $\mathcal{E}(a)$.

Consider, for example, the case of a hypersurface $X \subset M$. Let $\text{inv}_{1/2}(a) = \nu_1(a)$ denote the order $\nu_{X,a}$ of $X$ at a point $a$. Suppose that $g(x) = 0$ is a local defining equa-tion of $X$ at $a$ (i.e., $g$ generates $\mathcal{I}_{X,a}$). Let $N(a) =$ the germ of $M$ at $a$, $\mathcal{G}(a) = \{(g, d)\}$, where $d = \nu_1(a)$, and $\mathcal{E}(a) = \emptyset$. Then $\big(N(a), \mathcal{G}(a) = \mathcal{G}_1(a), \mathcal{E}(a)\big)$ is a codimension zero presentation of $\nu_1$ at $a$. Our local construction provides a way to define $\nu_2$ (and, in general, the successive $\nu_{r+1}$) by induction on codimension; the key point is that we can choose $z \in \mathcal{O}_{M,a}$ such that $\mu_a(z) = 1$ and $(N(a), \mathcal{G}(a), \mathcal{E}(a))$ is equivalent to $(N(a), \mathcal{G}(a) \cup \{(z, 1)\}, \mathcal{E}(a))$. It follows that, after any sequence of transformations of types (i), (ii) and (iii), $\mu_{a'}(z') = 1$, $S_{\mathcal{G}(a')} \subset \{z' = 0\}$ and $(N(a'), \mathcal{G}(a'), \mathcal{E}(a'))$ is equivalent to $(N(a'), \mathcal{G}(a') \cup \{(z', 1)\}, \mathcal{E}(a'))$ (Proposition 4.12).

The element $z$ above can be constructed explicitly as follows: Suppose that $(x_1, \ldots, x_n)$ is a local coordinate system for $M$ at $a$. (After a linear change of coordinates) we can assume that $(\partial^d g/\partial x_n^d)(a) \neq 0$. Then we can take $z = \partial^{d-1} g/\partial x_n^{d-1}$. By the implicit function theorem, $z = 0$ defines a (germ at $a$ of) a regular submanifold $N_1 = N_1(a)$ of $M$ of codimension 1. If $\mathcal{C}_1(a)$ denotes the collection of pairs $(h, \mu_h) = \left(\dfrac{\partial^q g}{\partial x_n^q}\Big|N_1,\, d - q\right)$, $q = 0, \ldots, d-2$ (each $h$ makes sense as an element of $\mathcal{O}_{N_1,a}$), then $\big(N_1(a), \mathcal{C}_1(a), \mathcal{E}(a) = \emptyset\big)$ is an (equivalent) codimension 1 presentation of $\nu_1$ at $a$.

Now consider any sequence (1.1) of blowings-up $\sigma_{j+1}$ with $\frac{1}{2}$-admissible centres $C_j$. Suppose that $a \in M_j$. Let $i$ be the smallest index $k$ such that $\nu_1(a) = \nu_1(a_k)$, as before; in particular, $E(a_i) = E^1(a_i)$. Let $(N(a_i), \mathcal{G}_1(a_i), \mathcal{E}(a_i) = \emptyset)$ be a codimension zero presentation of $\nu_1$ at $a_i$, and let $(N(a), \mathcal{G}_1(a), \mathcal{E}(a))$ denote its transform at $a$ (by the sequence of admissible blowings-up $\sigma_{i+1}, \ldots, \sigma_j$). Then $\mathcal{E}(a) = E(a) \backslash E^1(a)$ ($= \mathcal{E}_1(a)$, say), and $(N(a), \mathcal{G}_1(a), \mathcal{E}_1(a))$ is a codimension zero presentation of $\nu_1$ at $a$ with respect to $\mathcal{E}_1(a)$. For each $H \in E^1(a)$, let $\ell_H \in \mathcal{O}_{M_j,a}$ denote a generator of $\mathcal{I}_{H,a}$, and let $\mathcal{F}_1(a)$ denote $\mathcal{G}_1(a)$ together with all pairs $(f, \mu_f) = (\ell_H, 1)$, $H \in E^1(a)$. Then $(N(a), \mathcal{F}_1(a), \mathcal{E}_1(a))$ is a codimension zero presentation of $\mathrm{inv}_1 = (\nu_1, s_1)$ at $a$.

As above, choose $z \in \mathcal{O}_{M_i,a_i}$ such that $\mu_{a_i}(z) = 1$ and $(N(a_i), \mathcal{G}_1(a_i), \mathcal{E}_1(a_i))$ is equivalent to $(N(a_i), \mathcal{G}_1(a_i) \cup \{(z, 1)\}, \mathcal{E}_1(a_i))$. If $z'$ denotes the transform of $z$ at $a$, then $(N(a), \mathcal{G}_1(a), \mathcal{E}_1(a))$ is equivalent to $(N(a), \mathcal{G}_1(a) \cup \{(z', 1)\}, \mathcal{E}_1(a))$, and therefore $(N(a), \mathcal{F}_1(a), \mathcal{E}_1(a))$ is equivalent to $(N(a), \mathcal{F}_1(a) \cup \{(z', 1)\}, \mathcal{E}_1(a))$. Suppose that $(x_1, \ldots, x_n)$ is a local coordinate system for $M_j$ at $a$ such that $(\partial z'/\partial x_n)(a) \neq 0$. Let $N_1 = N_1(a)$ denote the (germ at $a$ of a) regular submanifold $\{z' = 0\}$, and let $\mathcal{H}_1(a)$ denote the collection of pairs $(h, \mu_h) = \left(\dfrac{\partial^q f}{\partial x_n^q}|N_1,\, \mu_f - q\right)$, $0 \leq q < \mu_f$, for all $(f, \mu_f) \in \mathcal{F}_1(a)$. Then $(N_1(a), \mathcal{H}_1(a), \mathcal{E}_1(a))$ is a codimension 1 presentation of $\mathrm{inv}_1$ at $a$. (Likewise, if $\mathcal{C}_1(a)$ denotes the collection of pairs $\left(\dfrac{\partial^q g}{\partial x_n^q}|N_1,\, \mu_g - q\right)$, $0 \leq q < \mu_g$, for all $(g, \mu_g) \in \mathcal{G}_1(a)$, then $(N_1(a), \mathcal{C}_1(a), \mathcal{E}_1(a))$ is a codimension 1 presentation of $\nu_1$ at $a$.)

Suppose that $\big(N_1(a), \mathcal{H}_1(a), \mathcal{E}_1(a)\big)$ is any codimension 1 presentation of $\mathrm{inv}_1$ at $a$, with respect to $\mathcal{E}_1(a) = E(a) \backslash E^1(a)$. Let $\mu_2(a) = \mu_{\mathcal{H}_1(a)}$. If $\mu_2(a) = \infty$, we set $\mathrm{inv}_X(a) = \big(\mathrm{inv}_1(a); \infty\big)$. Otherwise, for all $H \in \mathcal{E}_1(a)$, we write

$$\mu_{2H}(a) := \min\left\{\frac{\mu_{H,a}(h)}{\mu_h} \,:\, (h, \mu_h) \in \mathcal{H}_1(a)\right\},$$

where $\mu_{H,a}(h)$ denotes the *order of $h$ along $H \cap N_1$ at $a$* (i.e., the order to which a generator $x_H$ of the local ideal of $H \cap N_1$ factors from $h$); we define $\nu_2(a)$ as

$$\nu_2(a) := \mu_2(a) - \sum_H \mu_{2H}(a).$$

Then $\nu_2(a) \geq 0$. We will prove that each $\mu_{2H}(a)$ and thus $\nu_2(a)$ is an invariant of the equivalence class of $(N_1(a), \mathcal{H}_1(a), \mathcal{E}_1(a))$ (with respect to transformations (i), (ii) and (iii), but with a certain restriction on the sequence of transformations allowed; see Definition 4.10); hence each $\mu_{2H}(a)$ and $\nu_2(a)$ are invariants of the local isomorphism class of $M_j$, $X_j$, $E(a_j)$, $E^1(a_j)$.

Let $D(a) = \prod_{H \in \mathcal{E}(a)} x_H^{\mu_{2H}(a)}$, $D(a) = D_2(a)$; then each $h \in \mathcal{H}_1(a)$ can be factored as $h = D^{\mu_h} \cdot g$, and $\mu_a(g) \geq \mu_g$, where $\mu_g = \mu_h \cdot \nu_2(a)$. (Rational exponents and rational orders can be avoided by raising to suitable powers.) Let $\mathcal{G}_2(a)$ denote the collection of pairs $\{(g, \mu_g)\}$ together with $(D, 1 - \nu_2(a))$ in the case that $\nu_2(a) < 1$. ($\mathcal{G}_2(a) := \{(D, 1)\}$ in the case that $\nu_2(a) = 0$.) It follows that $(N_1(a), \mathcal{G}_2(a), \mathcal{E}_1(a))$ is a codimension 1 presentation of $\mathrm{inv}_{1\frac{1}{2}}$ at $a$ with respect to $\mathcal{E}_1(a) = E(a) \backslash E^1(a)$. If $\nu_2(a) = 0$, we set $\mathrm{inv}_X(a) = \mathrm{inv}_{1\frac{1}{2}}(a)$.

Suppose that $0 < \nu_2(a) < \infty$. Clearly, $\mu_{\mathcal{G}_2(a)} = 1$. Now assume that the blowings-up $\sigma_{j+1}$ in (1.1) are $1\frac{1}{2}$-admissible. Set $\mathcal{E}_2(a) = \mathcal{E}_1(a) \backslash E^2(a)$. Then $(N_1(a), \mathcal{G}_2(a), \mathcal{E}_2(a))$ is a codimension 1 presentation of $\mathrm{inv}_{1\frac{1}{2}}$ at $a$ with respect to $\mathcal{E}_2(a)$ and, as above, there is an equivalent codimension 2 presentation $(N_2(a), \mathcal{C}_2(a), \mathcal{E}_2(a)), \ldots$. The construction above can be repeated in increasing codimension. Eventually we reach $t \leq n = \dim_a M_j$ such that $0 < \nu_r(a) < \infty$ if $r \leq t$, and $\nu_{t+1}(a) = 0$ or $\infty$. Then we define $\mathrm{inv}_X(a) = (\mathrm{inv}_t(a); \nu_{t+1}(a))$ and $\mu_X(a) = \mu_{t+1}(a)$. See Chapter II. Our presentations satisfy a natural property of "semicoherence" (6.4) which allows us to prove that $\mathrm{inv}_X$ is Zariski-semicontinuous using the (elementary) Zariski-semicontinuity of order of a regular function. In Chapter II, we thus prove Theorem 1.14 (and therefore resolution of singularities) in the case of a hypersurface.

*Remark 1.18.* In the context of Remark 1.8, we can obtain a codimension zero presentation $(N(a), \mathcal{G}(a), \mathcal{E}(a) = \emptyset)$ of $\nu_1 = \nu_{\mathcal{I}}$ at $a$ (with respect to the notion of weak transform) simply by taking $N(a) =$ the germ of $M$ at $a$, and $\mathcal{G}(a) = \{(g, \nu_1(a))\}$, where $\{g\}$ is any finite set of generators of $\mathcal{I}_a$. The construction above allows us to define $\mathrm{inv}_{\mathcal{I}}(\cdot)$ and thus to prove the analogue of Theorem 1.14, and Theorem 1.10.

*Presentation of the Hilbert-Samuel function*

In higher codimension, the local construction described above can be used to define $\mathrm{inv}_X$ exactly as in the case of a hypersurface, provided that we can find a (semicoherent) presentation of the Hilbert-Samuel function. This is the problem solved in Chapter III. The standard basis of $\widehat{\mathcal{I}}_{X,a} \subset \widehat{\mathcal{O}}_{M,a}$ (with respect to any identification $\widehat{\mathcal{O}}_{M,a} \cong \underline{k}[[X_1, \ldots, X_n]]$) provides a *formal* presentation of $H_{X,\cdot}$ at $a$. The Henselian division theorem of Hironaka [H3] provides a presentation (at least with respect to admissible blowings-up (i); cf. [BM4, Th. 7.3]) which is algebraic in the sense of Artin, and therefore involves passing to an étale covering of $X$. We use a completely elementary division algorithm to obtain a presentation by regular functions. We also give $S_{H_{X,a}} = \{x \in |X| : H_{X,x} \geq H_{X,a}\}$ a natural structure of a closed *subspace* of $X$, and prove the coincidence of the sheaves of ideals defining $S_{H_{X,a}}$, on the one hand, and the equimultiple locus of a regular presentation, on the other.

*Remark 1.19.* The standard basis of $\widehat{\mathcal{I}}_{X,a}$ itself extends to a presentation of the Hilbert-Samuel function which is regular in a weaker sense that nevertheless suffices to prove desingularization using Chapter II: Let $\{F\} \subset \widehat{\mathcal{I}}_{X,a}$ denote the standard basis (with respect to a generic coordinate system). Then all formal derivatives $\partial^{|\alpha|} F/\partial X^\alpha$, $\alpha \in \mathbb{N}^n$, when restricted to $S_{H_{X,a}}$, are (induced by) regular functions defined in a common neighbourhood $U$ of $a$. Moreover, the induced formal expansions at each $b \in S_{H_{X,a}} \cap U$ provide a formal presentation at $b$. (This was our original approach and seems of independent interest; we plan to publish details elsewhere.)

The construction of our presentation of the Hilbert-Samuel function is based on combinatorial properties of the "diagram of initial exponents". (See §3 and (7.1) below.) Some of these techniques were developed in our earlier papers on resolution of singularities [BM4] and on differentiable functions and subanalytic sets [BM1], [BM2]. The stabilization theorem of §8 below is a new combinatorial result which plays an important part.

## 2. Examples and an application

In the examples below, we will follow the desingularization algorithm (over a sequence of local blowings-up of a hypersurface) as sketched in "Presentation of the invariant" in §1 and detailed in Chapter II. We will use the notation from §1.

*Example 2.1.* Consider the hypersurface $X = V(g) \subset \underline{k}^3$ defined by $g(x) = x_3^2 - x_1^2 x_2^3$.

*Year zero.* Let $a = 0$. Then $\nu_1(a) = \mu_a(g) = 2$ and $E(a) = \emptyset$, so that $s_1(a) = 0$. A codimension zero presentation of $\text{inv}_{1/2} = \nu_1$ at $a$ (with respect to $\mathcal{E}_1(a) = \emptyset$) is given by $(N(a), \mathcal{G}_1(a), \mathcal{E}_1(a) = \emptyset)$, where $N(a) = \underline{k}^3$ and $\mathcal{G}_1(a) = \{(g, 2)\} = \mathcal{F}_1(a)$. We can take $N_1(a) = \{x_3 = 0\}$ and $\mathcal{H}_1(a) = \{(x_1^2 x_2^3, 2)\}$ to get a codimension 1 presentation $(N_1(a), \mathcal{H}_1(a), \mathcal{E}_1(a))$ of $\text{inv}_1 = (\nu_1, s_1)$ at $a$. Thus, $\nu_2(a) = \mu_2(a) = 5/2$ and $\text{inv}_{1\frac{1}{2}}(a) = (2, 0; 5/2)$. Let $\mathcal{G}_2(a) = \{(x_1^2 x_2^3, 5)\}$; then $(N_1(a), \mathcal{G}_2(a), \mathcal{E}_1(a))$ is a codimension 1 presentation of $\text{inv}_{1\frac{1}{2}}$ at $a$. The latter presentation is equivalent to $(N_1(a), \{(x_1, 1), (x_2, 1)\}, \mathcal{E}_1(a) = \emptyset)$, so repeating the construction, we compute $\text{inv}_X(a) = (2, 0; 5/2, 0; 1, 0; \infty)$ and $S_{\text{inv}_X}(a) = S_{\text{inv}_{1\frac{1}{2}}}(a) = \{a\}$. ($S_{\text{inv}_X}(a)$ denotes the germ of $S_{\text{inv}_X(a)}$ at $a$, etc.) We therefore let $\sigma_1 : M_1 \to M_0 = \underline{k}^3$ be the blowing-up with centre $C_0 = \{a\}$. $M_1$ is covered by 3 coordinate charts $U_i = M_1 \setminus \{x_i = 0\}'$, where $\{x_i = 0\}'$ means the strict transform of $\{x_i = 0\}$, $i = 1, 2, 3$; $\sigma_1|U_1$ can be written $x_1 = y_1$, $x_2 = y_1 y_2$, $x_3 = y_1 y_3$ (cf. "Blowing up" in §3 below).

*Year one.* Let $X_1$ denote the strict transform of $X$ by $\sigma_1$; then $X_1 \cap U_1 = V(g_1)$, where $g_1 = g' := y_1^{-2} g \circ \sigma_1 = y_3^2 - y_1^3 y_2^3$. Let $b = 0$ in $U_1$. Then $\nu_1(b) = 2 = \nu_1(a)$; therefore, $E^1(b) = \emptyset$, $s_1(b) = 0$, and $\mathcal{E}_1(b) := E(b) \setminus E^1(b) = E(b) = \{H_1\}$, where $H_1 = \sigma_1^{-1}(a) = \{y_1 = 0\}$. We can take $\mathcal{F}_1(b) = \mathcal{G}_1(b) = \{(g_1, 2)\}$, $N_1(b) = \{y_3 = 0\} = N_1(a)'$, and $\mathcal{H}_1(b) = \{(y_1^3 y_2^3, 2)\}$. (Of course $y_1^3 y_2^3 = y_1^{-2}((x_1^2 x_2^3) \circ \sigma_1)$.) Then $\mu_2(b) = 3$ and $\mu_{2H_1}(b) = 3/2$, so that $\nu_2(b) = 3 - 3/2 = 3/2$ and $\text{inv}_{1\frac{1}{2}}(b) = (2, 0; 3/2)$. $D_2(b) = y_1^{3/2}$, so that $\mathcal{G}_2(b) = \{(y_2^3, 3)\}$, which is equivalent to $\{(y_2, 1)\}$. $(N_1(b), \mathcal{G}_2(b), \mathcal{E}_1(b))$ is a presentation of $\text{inv}_{1\frac{1}{2}}$ at $b$; therefore, $S_{\text{inv}_{1\frac{1}{2}}}(b) = \{y_2 = y_3 = 0\}$. Repeating the procedure: $E^2(b) = \{H_1\}$, $\text{inv}_2(b) = (2, 0; 3/2, 1)$ and $\text{inv}_2$ is presented at $b$ by $(N_1(b), \mathcal{F}_2(b), \mathcal{E}_2(b) = \emptyset)$, where $\mathcal{F}_2(b) = \{(y_1, 1), (y_2, 1)\}$. Finally, $\text{inv}_X(b) = (2, 0; 3/2, 1; 1, 0; \infty)$ and $S_{\text{inv}_X}(b) = $ $S_{\mathrm{inv}_2}(b) = \{y_1 = y_2 = y_3 = 0\} = \{b\}$. We let $\sigma_2$ be the blowing-up with centre $C_1 = \{b\}$. $\sigma_2^{-1}(U_1)$ is covered by 3 coordinate charts $U_{1i} = \sigma_2^{-1}(U_1)\backslash\{y_i = 0\}'$, $i = 1, 2, 3$; $\sigma_2|U_{12}$ can be written $y_1 = z_1 z_2$, $y_2 = z_2$, $y_3 = z_2 z_3$.

*Year two.* Let $X_2$ denote the strict transform of $X_1$; in particular, $X_2 \cap U_{12} = V(g_2)$, where $g_2 = g_1' = z_2^{-2} g_1 \circ \sigma_2 = z_3^2 - z_1^3 z_2^4$. Let $c = 0$ in $U_{12}$. Now, $E(c) = \{H_1, H_2\}$, where $H_1 = \{y_1 = 0\}' = \{z_1 = 0\}$ and $H_2 = \sigma_2^{-1}(b) = \{z_2 = 0\}$. Then $\nu_1(c) = 2 = \nu_1(a)$, so that $E^1(c) = \emptyset$, $s_1(c) = 0$, and $\mathcal{E}_1(c) = E(c)$. We take $\mathcal{F}_1(c) = \mathcal{G}_1(c) = \{(g_2, 2)\}$, $N_1(c) = \{z_3 = 0\}$, and $\mathcal{H}_1(c) = \{(z_1^3 z_2^4, 2)\}$. (Again, $z_1^3 z_2^4 = z_2^{-2}((y_1^3 y_2^3) \circ \sigma_2)$.) Then $\mu_2(c) = 7/2$ and $D_2(c) = z_1^{3/2} z_2^2$, so that $\nu_2(c) = 0$ and $\mathrm{inv}_X(c) = \mathrm{inv}_{1\frac{1}{2}}(c) = (2, 0; 0)$. $(N_1(c), \mathcal{H}_1(c), \mathcal{E}_1(c))$ is a presentation of $\mathrm{inv}_1$ (or of $\mathrm{inv}_{1/2}$) at $c$, and $(N_1(c), \mathcal{G}_2(c), \mathcal{E}_1(c))$, where $\mathcal{G}_2(c) = \{(D_2(c), 1)\}$, is a presentation of $\mathrm{inv}_X = \mathrm{inv}_{1\frac{1}{2}}$ at $c$. $S_{\mathrm{inv}_X}(c) = S_{\mathrm{inv}_{1\frac{1}{2}}}(c)$ is the union of the $z_2$- and $z_1$-axes; $S_{\mathrm{inv}_X}(c) \cap H_1 = z_2$-axis and $S_{\mathrm{inv}_X}(c) \cap H_2 = z_1$-axis. In the lexicographic ordering of the set of subsets of $E(c)$ (given by Remark 1.15 (2)), $\{H_1\} = (1, 0) > (0, 1) = \{H_2\}$, so that $J(c) = \{H_1\}$ and $\mathrm{inv}_X^e(c) = (\mathrm{inv}_X(c); \{H_1\})$. In other words, although (for property (4) of Theorem 1.14) we could choose either component of $S_{\mathrm{inv}_X}(c)$ as the next centre of blowing-up, for the purpose of canonical desingularization we choose $C_2 = z_2$-axis. Let $\sigma_3$ by the blowing-up with centre $C_2$. $\sigma_3^{-1}(U_{12})$ is covered by 2 coordinate charts $U_{12i}$, where $U_{12i} = \sigma_3^{-1}(U_{12})\backslash\{z_i = 0\}'$, $i = 1, 3$; $\sigma_3|U_{121}$ can be written $z_1 = w_1$, $z_2 = w_2$, $z_3 = w_1 w_3$.

*Year three.* Let $X_3$ be the strict transform of $X_2$; in particular, $X_3 \cap U_{121} = V(g_3)$, where $g_3 = g_2' = w_3^2 - w_1 w_2^4$. Let $d = 0$ in $U_{121}$, so that $E(d) = \{H_2, H_3\}$, where $H_2 = \{z_2 = 0\}' = \{w_2 = 0\}$ and $H_3 = \sigma_3^{-1}(C_2) = \{w_1 = 0\}$. Then $\nu_1(d) = 2$, so that $E^1(d) = \emptyset$, $s_1(d) = 0$, and $\mathcal{E}_1(d) = E(d)$. We take $N_1(d) = \{w_3 = 0\}$ (still the strict transform of $N_1(a) = \{x_3 = 0\}$) and $\mathcal{H}_1(d) = \{(w_1 w_2^4, 2)\}$. Then $\mu_2(d) = 5/2$ and $D_2(d) = w_1^{1/2} w_2^2$, so that $\nu_2(d) = 0$ and $\mathrm{inv}_X(d) = \mathrm{inv}_{1\frac{1}{2}}(d) = (2, 0; 0)$. Again, $(N_1(d), \mathcal{H}_1(d), \mathcal{E}_1(d))$ is a presentation of $\mathrm{inv}_1$ (or of $\mathrm{inv}_{1/2}$) at $d$, and $(N_1(d), \mathcal{G}_2(d), \mathcal{E}_1(d))$, where $\mathcal{G}_2(d) = \{(D_2(d), 1)\}$ is a presentation of $\mathrm{inv}_X = \mathrm{inv}_{1\frac{1}{2}}$ at $d$. Therefore, $S_{\mathrm{inv}_X}(d) = S_{\mathrm{inv}_{1\frac{1}{2}}}(d) = \{w_2 = w_3 = 0\}$. We let $\sigma_4$ be the blowing-up with centre $C_3 = w_1$-axis. Note that $\mathrm{inv}_X(d) = \mathrm{inv}_X(c)$, but $\mu_X(d) = \mu_2(d) = 5/2 < 7/2 = \mu_2(c) = \mu_X(c)$ (as predicted by Theorem 1.14 (4)). $\sigma_4^{-1}(U_{121})$ is covered by 2 coordinate charts $U_{121i}$, where $U_{121i} = \sigma_4^{-1}(U_{121})\backslash\{w_i = 0\}'$, $i = 2, 3$; $\sigma_4|U_{1212}$ can be written $w_1 = v_1$, $w_2 = v_2$, $w_3 = v_2 v_3$.

*Year four.* Let $X_4$ be the strict transform of $X_3$; thus $X_4 \cap U_{1212} = V(g_4)$, where $g_4 = g_3' = v_3^2 - v_1 v_2^2$. Let $e = 0$ in $U_{1212}$, so that $E(e) = \{H_3, H_4\}$, where $H_3 = \{w_1 = 0\}' = \{v_1 = 0\}$ and $H_4 = \sigma_4^{-1}(C_3) = \{v_2 = 0\}$. Then $\nu_1(e) = 2$, so that $E^1(e) = \emptyset$, $s_1(e) = 0$, $\mathcal{E}_1(e) = E(e)$. As above, we get $\mu_2(e) = 3/2$ and $D_2(e) = v_1^{1/2} v_2$, so that $\mathrm{inv}_X(e) = (2, 0; 0)$. $\mathrm{inv}_X$ is presented at $e$ by $(N_1(e), \mathcal{G}_2(e), \mathcal{E}_1(e))$, where $N_1(e) = \{v_3 = 0\}$ and $\mathcal{G}_2(e) = \{(D_2(e), 1)\}$. Therefore, $S_{\mathrm{inv}_X}(e) = \{v_2 = v_3 = 0\}$. It is easy to see that, if we blow up with centre $C_4 = S_{\mathrm{inv}_X}(e)$, then the multiplicity of the strict transform deceases; in fact, the strict transform $X_5$ is non-singular.

*Example 2.2.* Consider $X = \{x_3^2 - x_1 x_2^2 = 0\} \sim$ the same hypersurface as in year four above $\sim$ but without a history of previous blowings-up; i.e., $E(\cdot) = \emptyset$ everywhere. Let $a = 0$. In this case, $\mathrm{inv}_{1\frac{1}{2}}(a) = (2, 0; 3/2)$ (cf. year zero above), and we can take $N_1(a) = \{x_3 = 0\}$, $\mathcal{H}_1(a) = \{(x_1 x_2^2, 0)\}$ and $\mathcal{G}_2(a) = \{(x_1 x_2^2, 3)\}$; $(N_1(a), \mathcal{G}_2(a), \mathcal{E}_1(a) = \emptyset)$ is a codimension 1 presentation of $\text{inv}_{1\frac{1}{2}}$ at $a$, and we get an equivalent presentation by replacing $\mathcal{G}_2(a)$ with $\{(x_1, 1), (x_2, 1)\}$. Therefore, $\text{inv}_X(a) = (2, 0; 3/2, 0; 1, 0; \infty)$ (as in year zero above) and $S_{\text{inv}_X}(a) = S_{\text{inv}_{1\frac{1}{2}}}(a) = \{a\}$. As centre of blowing up we would choose $C = S_{\text{inv}_X}(a) = \{a\} \sim$ not the $x_1$-axis as in year four of Example 2.1, although the singularity is the same!

*Example 2.3.* Consider the hypersurface $X = V(g) \subset \underline{k}^3$, where $g(x) = x_3^3 - x_1 x_2$.

*Year zero.* Let $a = 0$. Then $\nu_1(a) = \mu_a(g) = 2$ and $\text{Sing} X = \{a\}$, so that $S_{\text{inv}_X}(a) = \{a\}$. We therefore let $\sigma_1 : M_1 \to M_0 = \underline{k}^3$ be the blowing-up with centre $C_0 = \{a\}$. $M_1$ is covered by 3 coordinate charts $U_i = M_1 \setminus \{x_i = 0\}'$, where $\{x_i = 0\}'$ is the strict transform of $\{x_i = 0\}$, $i = 1, 2, 3$; $\sigma_1 | U_3$ can be written $x_1 = y_1 y_3$, $x_2 = y_2 y_3$, $x_3 = y_3$.

*Year one.* Let $X_1$ denote the strict transform of $X$ by $\sigma_1$; then $X_1 \cap U_3 = V(g_1)$, where $g_1 = y_3^{-2} g \circ \sigma_1 = y_3 - y_1 y_2$. Let $b = 0$ in $U_3$. Then $\nu_1(b) = 1 < 2 = \nu_1(a)$; therefore $E^1(b) = E(b) = \{H_1\}$, where $H_1 = \sigma_1^{-1}(a) = \{y_3 = 0\}$, so that $s_1(b) = 1$ and $\mathcal{E}_1(b) = \emptyset$. We can take $\mathcal{F}_1(b) = \{(g_1, 1), (y_3, 1)\}$, $N_1(b) = \{y_3 = 0\}$ and $\mathcal{H}_1(b) = \{(y_1 y_2, 1)\}$. Then $\mu_2(b) = 2 = \nu_2(b)$, $\text{inv}_{1\frac{1}{2}}(b) = (1, 1; 2)$ and $\mathcal{F}_2(b) = \mathcal{G}_2(b) = \{y_1 y_2, 2)\}$, which is equivalent to $\{(y_1, 1), (y_2, 1)\}$. It follows that $\text{inv}_X(b) = (1, 1; 2, 0; 1, 0; \infty)$ and $S_{\text{inv}_X}(b) = S_{\text{inv}_{1\frac{1}{2}}}(b) = \{b\}$. Let $\sigma_2$ be the blowing-up with centre $C_1 = \{b\}$. $\sigma_2^{-1}(U_3)$ is covered by 3 coordinate charts $U_{3i} = \sigma_2^{-1}(U_3) \setminus \{y_i = 0\}'$, $i = 1, 2, 3$; $\sigma_2 | U_{31}$ can be written $y_1 = z_1$, $y_2 = z_1 z_2$, $y_3 = z_1 z_3$.

*Year two.* Let $X_2$ be the strict transform of $X_1$; in particular, $X_2 \cap U_{31} = V(g_2)$, where $g_2 = z_1^{-1} g_1 \circ \sigma_2 = z_3 - z_1 z_2$. Let $c = 0$ in $U_{31}$. Then $\nu_1(c) = 1 = \nu_1(b)$, and $E(c) = \{H_1, H_2\}$, where $H_1 = \{y_3 = 0\}' = \{z_3 = 0\}$ and $H_2 = \sigma_2^{-1}(b) = \{z_1 = 0\}$, so that $E^1(c) = \{H_1\}$, $s_1(c) = 1$ and $\mathcal{E}_1(c) = \{H_2\}$. We take $\mathcal{F}_1(c) = \{(g_2, 1), (z_3, 1)\}$, $N_1(c) = \{z_3 = 0\}$ and $\mathcal{H}_1(c) = \{(z_1 z_2, 1)\}$. Then $\mu_2(c) = 2$ and $D_2(c) = z_1$, so that $\nu_2(c) = 1$ and $\text{inv}_{1\frac{1}{2}}(c) = (1, 1; 1)$. Hence $E^2(c) = \{H_2\}$ and $(N_1(c), \mathcal{G}_2(c), \mathcal{E}_2(c))$, where $\mathcal{G}_2(c) = \{(z_2, 1)\}$ and $\mathcal{E}_2(c) = \emptyset$, is a presentation of $\text{inv}_{1\frac{1}{2}}$ at $c$. It follows that $\text{inv}_X(c) = (1, 1; 1, 1; 1, 0; \infty)$ and $S_{\text{inv}_X}(c) = S_{\text{inv}_2}(c) = \{c\}$. We therefore let $\sigma_3$ be the blowing-up with centre $C_2 = \{c\}$. $\sigma_3^{-1}(U_{31})$ is covered by 3 coordinate charts $U_{31i} = \sigma_3^{-1}(U_{31}) \setminus \{z_i = 0\}'$, $i = 1, 2, 3$; $\sigma_3 | U_{311}$ can be written $z_1 = w_1$, $z_2 = w_1 w_2$, $z_3 = w_1 w_3$.

*Year three.* Let $X_3$ be the strict transform of $X_2$; in particular, $X_3 \cap U_{311} = V(g_3)$, where $g_3 = w_1^{-1} g_2 \circ \sigma_3 = w_3 - w_1 w_2$. Let $d = 0$ in $U_{311}$, so that $E(d) = \{H_1, H_3\}$, where $H_1 = \{w_3 = 0\}$ and $H_3 = \sigma_3^{-1}(c) = \{w_1 = 0\}$. Then $\nu_1(d) = 1$, $E^1(d) = \{H_1\}$, $s_1(d) = 1$ and $\mathcal{E}_1(d) = \{H_3\}$. We take $\mathcal{F}_1(d) = \{(g_3, 1), (w_3, 1)\}$, $N_1(d) = \{w_3 = 0\}$ and $\mathcal{H}_1(d) = \{(w_1 w_2, 1)\}$. Then $\mu_2(d) = 2$ and $D_2(d) = w_1$, so that $\nu_2(d) = 1$ and $\text{inv}_{1\frac{1}{2}}(d) = (1, 1; 1)$. Hence $E^2(d) = \emptyset$, $\text{inv}_2(d) = (1, 1; 1, 0)$ and $(N_1(d), \mathcal{F}_2(d), \mathcal{E}_2(d))$, where $\mathcal{F}_2(d) = \mathcal{G}_2(d) = \{(w_2, 1)\}$ and $\mathcal{E}_2(d) = \{H_3\}$, is a presentation of $\text{inv}_2$ at $d$. It follows that $\text{inv}_X(d) = (1, 1; 1, 0; \infty)$ and $S_{\text{inv}_X}(d) = \{w_3 = 0, w_2 = 0\}$. In this chart $U_{311}$, $X_3$ is smooth and has only normal crossings simultaneously with respect to the collection $E_3$ of all exceptional divisors at every point of $\{w_3 = w_2 = 0\}$ except $d = 0$ (cf. Remarks 1.7(3)).

*An application: Łojasiewicz's inequalities*

The fundamental inequalities of Łojasiewicz are immediate consequences of desingularization in the form of Theorem 1.10 (or Theorem 1.6 in the hypersurface case); in fact, we need only the following:

**Theorem 2.4.** *Let $M$ be a manifold, and let $\mathcal{I} \subset \mathcal{O}_M$ denote a sheaf of (principal) ideals of finite type. Then there is a manifold $M'$ and a proper surjective morphism $\varphi : M' \to M$ such that $\varphi^{-1}(\mathcal{I})$ is a normal-crossings divisor.*

**Theorem 2.5. Inequality I.** *Assume that $\underline{k} = \mathbb{R}$ or $\mathbb{C}$. Let $f$ and $g$ denote regular functions on a manifold $M$. (Recall that "regular" means "analytic" in the category of analytic spaces.) Suppose that $K$ is a compact subset of $M$ and that $\{x : g(x) = 0\} \subset \{x : f(x) = 0\}$ in a neighbourhood of $K$. Then there are positive constants $c$, $\lambda$ such that*

$$|g(x)| \geq c|f(x)|^\lambda$$

*in some neighbourhood of $K$. Moreover, the infimum of such $\lambda$ is a positive rational number.*

*Proof.* This is obvious if $f(x) \cdot g(x)$ has only normal crossings in a neighbourhood of $K$; in general, therefore, it follows from Theorem 2.4. □

*Remark 2.6.* We are assuming here that the category of spaces is from (0.2) (2) or (3). (If $M$ has a quasi-compact underlying algebraic structure with respect to which $f$ and $g$ are regular, then $\lambda$ can be chosen independent of $K$; there is an analogous remark concerning Inequalities II and III following.) The argument above allows us to conclude that, in any of the categories of (0.2), locally some power of $f \circ \varphi$ belongs to the ideal generated by $g \circ \varphi$; it follows that locally $f$ belongs to the integral closure of the ideal generated by $g$, and the equation of integral dependence has degree bounded on $K$ (cf. [LT]).

**Theorem 2.7. Inequality II.** *Let $f$ be a regular function on an open subspace $M$ of $\mathbb{R}^n$. Suppose that $K$ is a compact subset of $M$, on which $\mathrm{grad}\, f(x) = 0$ only if $f(x) = 0$. Then there exists $c > 0$ and $\mu$, $0 < \mu \leq 1$, such that*

$$|\mathrm{grad}\, f(x)| \geq c|f(x)|^{1-\mu}$$

*in a neighbourhood of $K$. (Sup$\mu$ is rational.)*

*Proof.* Note that, if $f(a) = 0$, then there is a neighbourhood of $a$ in which $\mathrm{grad} f(x) = 0$ only if $f(x) = 0$. Let $g(x) = |\mathrm{grad} f(x)|^2 = \sum_{i=1}^n (\partial f/\partial x_i)^2$. (As in the proof of Inequality I) let $\varphi : M' \to M$ be a morphism given by Theorem 2.4 for the ideal generated by $f \cdot g$. We claim there is a neighbourhood of $\varphi^{-1}(K)$ in which $\varphi^*(f^2/g)$ is a regular function vanishing on $\{x : (f \circ \varphi)(x) = 0\}$:

Consider any regular curve $\gamma : x = x(t)$ in $M$ such that $\gamma \cap \{x : f(x) = 0\} = \{x(0)\}$ and $\gamma$ is the image of a smooth curve in $M'$ which is transverse to $\varphi^{-1}(\{x : g(x) = 0\})$ at a smooth point of the latter. Let $Q(t) = f(x(t))$. Then $Q(t) \neq 0$ for $t \neq 0$, so that $Q(t)$ has nonzero Taylor expansion at $t = 0$. Therefore, $Q(t)$ is divisible by $Q'(t) = df(dx/dt)$ and the quotient $Q(t)/Q'(t)$ vanishes at $t = 0$. Since $|Q'(t)|^2 \leq g(x(t))|dx/dt|^2$, it follows that $Q(t)^2$ is divisible by $g(x(t))$ and the quotient $f(x(t))^2/g(x(t))$ vanishes at $t = 0$. The claim follows.

From the claim, we can conclude (as in Theorem 2.5) that there are positive constants $c$ and $\mu$ (where $\sup \mu$ is rational) such that

$$|f(x)|^{2\mu} \geq c^2 \frac{f(x)^2}{g(x)}$$

in a neighbourhood of $K$. Clearly, $0 < \mu \leq 1$ ($\mu = 1$ if and only if $g(x)$ vanishes nowhere on $K$.) Thus,

$$|\mathrm{grad} f(x)| \geq c|f(x)|^{1-\mu} \ . \qquad \Box$$

**Theorem 2.8. Inequality III.** *Let $f$ be a regular function on an open subspace $M$ of $\mathbb{R}^n$, and set $Z = \{x \in M : f(x) = 0\}$. Suppose that $K$ is a compact subset of $M$. Then there are $c > 0$ and $\nu \geq 1$ such that*

$$|f(x)| \geq c d(x, Z)^{\nu}$$

*in a neighbourhood of $K$. ($d(\cdot, Z)$ denotes the distance to $Z$.) Again the infimum of such $\nu$ is rational.*

*Proof.* This follows from Inequality II: We can assume that $\mathrm{grad} f(x) = 0$ only if $f(x) = 0$, on $K$. We then claim that (even if $f$ is merely $\mathcal{C}^1$ and) if $|\mathrm{grad} f(x)| \geq c|f(x)|^{1-\mu}$ in a neighbourhood $U$ of $K$, where $0 < \mu \leq 1$, then $|f(x)|^{\mu} \geq \mu c d(x, Z)$ in some neighbourhood of $K$ (cf. Łojasiewicz [Ł]):

Consider a point $a \in U$ such that $f(a) \neq 0$. We can assume that $f(a) > 0$. (Otherwise, use $-f$.) Suppose that $x(t)$ is a solution of the equation

$$\frac{dx}{dt} = -\frac{\mathrm{grad} f(x)}{|\mathrm{grad} f(x)|}$$

with $x(0) = a$. Write $Q(t) = f(x(t))$. Then $Q'(t) = df(dx/dt) = -|\mathrm{grad} f(x(t))| < 0$. Hence

$$\begin{aligned}
\frac{f(a)^{\mu}}{\mu} &\geq \frac{Q(0)^{\mu} - Q(t)^{\mu}}{\mu} \\
&= -\frac{1}{\mu} \int_0^t \frac{d}{dt} Q(t)^{\mu} dt \\
&= -\int_0^t Q(t)^{\mu-1} Q'(t) dt \\
&\geq c \int_0^t dt = ct \ .
\end{aligned}$$

It follows that the solution curve $x = x(t)$ tends to $Z$ in a finite time $t_0$. Since $|dx/dt| = 1$, $t_0 \geq d(a, Z)$ and $f(a)^{\mu} \geq \mu c d(a, Z)$, as required. $\qquad \Box$

## 3. Basic notions

*Definitions and notation*

Let $X = (|X|, \mathcal{O}_X)$ denote a local-ringed space over $\underline{k}$. We call $|X|$ the *support* or *underlying topological space* of $X$. We recall the following definitions: $X$ is *smooth* if, for all $x \in |X|$, $\mathcal{O}_{X,x}$ is a regular local ring. A local-ringed space $Y = (|Y|, \mathcal{O}_Y)$ is a *closed subspace* of $X$ if there is a sheaf of ideals $\mathcal{I}_Y$ of finite type in $\mathcal{O}_X$ such that $|Y| = \operatorname{supp} \mathcal{O}_X/\mathcal{I}_Y$ and $\mathcal{O}_Y$ is the restriction to $|Y|$ of $\mathcal{O}_X/\mathcal{I}_Y$. $Y$ is an *open subspace* of $X$ if $|Y|$ is an open subset of $|X|$ and $\mathcal{O}_Y = \mathcal{O}_X \big| |Y|$.

Let $X = (|X|, \mathcal{O}_X)$ be a local-ringed space. Let $a \in |X|$. Suppose that $f \in \mathcal{O}_{X,a}$ (or that $f \in \mathcal{O}_X(U)$, where $U$ is an open neighbourhood of $a$; we usually do not distinguish between $f \in \mathcal{O}_{X,a}$ and a representative in a suitable neighbourhood $U$). We define the *order* $\mu_a(f)$ *of* $f$ *at* $a$ as the largest $p \in \mathbb{N}$ such that $f \in \underline{m}_{X,a}^p$ (where $\underline{m}_{X,a}$ denotes the maximal ideal of $\mathcal{O}_{X,a}$). ($\mu_a f) = \infty$ if $f = 0$ in $\mathcal{O}_{X,a}$.) Let $C$ be a closed subspace of $X$, so that $C$ is defined by a sheaf of ideals $\mathcal{I}_C \subset \mathcal{O}_X$ of finite type. We define the *order* $\mu_{C,a}(f)$ *of* $f$ *along* $C$ *at* $a$ as the largest $p \in \mathbb{N}$ such that $f \in \mathcal{I}_{C,a}^p$.

Let $\varphi\colon X \to Y$ be a morphism of local-ringed spaces. Thus, for all $a \in |X|$, $\varphi$ induces a local homomorphism $\varphi_a^*\colon \mathcal{O}_{Y,\varphi(a)} \to \mathcal{O}_{X,a}$ (and a local homomorphism $\widehat{\varphi}_a^*\colon \widehat{\mathcal{O}}_{Y,\varphi(a)} \to \widehat{\mathcal{O}}_{X,a}$ of the completions). If $g \in \mathcal{O}_{Y,\varphi(a)}$ (or $g \in \mathcal{O}_Y(V)$, where $V$ is an open neighbourhood of $\varphi(a)$), then we will denote $\varphi_a^*(g)$ also by $g \circ \varphi_a$ or even by $g \circ \varphi$. (Similarly for $\widehat{\varphi}_a^*$.)

An element $f \in \mathcal{O}_X(U)$, where $U \subset |X|$ is open, will be called a *regular function* (*on* $U$). We will write $X_a$ (respectively, $f_a$) for the germ at $a$ of $X$ (respectively, of a regular function $f$). If $U$ is open in $|X|$ and $f_1, \ldots, f_\ell \in \mathcal{O}_X(U)$, then $V(f_1, \ldots, f_\ell)$ will denote the subspace of $X|U$ defined by the ideal subsheaf of $\mathcal{O}_X|U$ generated by the $f_i$.

*Regular coordinate charts*

If $M$ is an analytic manifold (i.e., $M = (|M|, \mathcal{O}_M)$ is a smooth analytic space over $\underline{k}$), then a classical coordinate chart $U$ is a regular coordinate chart in the sense of (0.3). (Here the ring of regular functions $\mathcal{O}(U) = \mathcal{O}_M(U)$ means the ring of analytic functions on $U$.)

In this subsection, we show how to construct regular coordinate charts in the algebraic context. Consider a scheme of finite type over $\underline{k}$. Let $X = (|X|, \mathcal{O}_X)$ denote either the scheme itself, or the local-ringed space where $|X|$ is the set of $\underline{k}$-rational points of the scheme, with the induced Zariski topology, and $\mathcal{O}_X$ is the restriction to $|X|$ of the structure sheaf of the scheme. We will show that if $X = M$ is smooth, then $M$ can be covered by coordinate charts as in (0.3). In the remainder of the article, we will adopt the convention that the residue field is $\underline{k}$ at every point (and we will write $\underline{k}^n$ rather than $\mathbb{A}^n$) in order to use a language that is common to schemes, analytic spaces, etc. But it will be clear from the construction of coordinate charts in this section, that all of our constructions apply to schemes of finite type over $\underline{k}$. We believe that our focus on the $\underline{k}$-rational points is very natural; for example, in resolving the singularities of real algebraic varieties.

Let $M = (|M|, \mathcal{O}_M)$. Each point $a$ of $M$ admits a Zariski-open neighbourhood $U$ in which regular functions (elements of $\mathcal{O}(U) = \mathcal{O}_M(U)$) can be described as follows:

(3.1) (1) $U = V(p_1, \ldots, p_{N-n})$, where $N \geq n = \dim_a M$ and the $p_j \in \underline{k}[u, v]$ are polynomials in $(u, v) = (u_1, \ldots, u_n, v_1, \ldots, v_{N-n})$ such that $\det \partial p/\partial v$ vanishes nowhere on $U$ (i.e., is invertible in the local ring of $\mathbb{A}^N$ at every point of $U$). ($\partial p/\partial v$ denotes the Jacobian matrix $\partial(p_1, \ldots, p_{N-n})/\partial(v_1, \ldots, v_{N-n})$.) We thus have a closed embedding $U \hookrightarrow \mathbb{A}^N$. (We say that the projection $(u, v) \mapsto u$ of $\mathbb{A}^N$ onto $\mathbb{A}^n$ induces an "étale covering" $U \to \mathbb{A}^n$.)

(2) Each element of $\mathcal{O}(U)$ is the restriction to $U$ of a rational function $f = q/r$, where $q, r \in \underline{k}[u, v]$ and $r$ vanishes nowhere on $U$.

In the case that $M$ is a scheme, $U = \operatorname{Spec} \underline{k}[u, v]/I$, where $I = (p_1, \ldots, p_{N-n})$ is the ideal generated by the $p_i$, and $\mathcal{O}(U)$ can be identified with $\underline{k}[u, v]/I$ (by the Nullstellensatz).

**Definition 3.2.** *A (regular) coordinate system $(x_1, \ldots, x_n)$ on a Zariski-open subset $U$ of $|M|$ means an $n$-tuple of elements $x_i \in \mathcal{O}(U)$ satisfying the following condition: Let $a \in U$. Let $a_i = x_i(a) \in \mathbb{F}_a$, $i = 1, \ldots, n$, where $\mathbb{F}_a$ denotes the residue field $\mathcal{O}_a/\underline{m}_a$. ($\underline{m}_a$ is the maximal ideal of $\mathcal{O}_a = \mathcal{O}_{M,a}$.) If $\Phi_i(z) \in \underline{k}[z]$ denotes the minimal polynomial of $a_i$ (i.e., the minimal monic relation for $a_i$ with coefficients in $\underline{k}$), $i = 1, \ldots, n$, then the $\Phi_i(x_i)$ form a basis of $\underline{m}_a/\underline{m}_a^2$ over $\mathbb{F}_a$.*

In this case, $\dim \mathcal{O}_a = \dim_{\mathbb{F}_a} \underline{m}_a/\underline{m}_a^2$. If $\mathbb{F}_a = \underline{k}$ (i.e., if $a$ is a $\underline{k}$-*rational point*) then $\Phi_i(x_i) = x_i - a_i$. In general, $\Phi_i(x_i) \sim x_i - a_i$ in the localization $\mathbb{F}_a[x_i]_{(a_i)}$. (We use $\sim$ to mean "= except for an invertible factor".)

In (3.1) above, for example, the restrictions $x_i$ to $U$ of the $u_i$ form a regular coordinate system $(x_1, \ldots, x_n)$. (The values of the coordinates may coincide at different points of $U$.)

**Lemma 3.3.** *Let $a \in M$ and let $x_1, \ldots, x_n$ denote regular functions on a neighbourhood of $a$. Then there is a Zariski-open neighbourhood $U$ of $a$ with the following property: $(x_1, \ldots, x_n)$ is a regular coordinate system on $U$ if and only if there is a closed embedding $U \hookrightarrow \mathbb{A}^N$ for some $N$, as in (3.1), such that the $x_i$ are the restrictions of the $u_i$ to $U$.*

*Proof.* Let $U$ be a Zariski-open neighbourhood of $a$ such that $U$ admits a closed embedding $U \hookrightarrow \mathbb{A}^N$ satisfying (3.1), and each $x_i \in \mathcal{O}(U)$; thus each $x_i$ is the restriction to $U$ of a rational function $q_i(u, v)/r_i(u, v)$, where $q_i, r_i \in \underline{k}[u, v]$ and $r_i(u, v)$ vanishes nowhere on $U$. Clearly, $(x_1, \ldots, x_n)$ forms a regular coordinate system on $U$ if and only if the gradients of the $q_i/r_i$ and the $p_j$ are linearly independent at every point of $U$. Consider

$$\begin{array}{ccc} U & \hookrightarrow & \mathbb{A}^{n+N} \quad (y, u, v) \\ & \searrow \downarrow & \downarrow \\ & \mathbb{A}^n & y \end{array}$$

where $y = (y_1, \ldots, y_n)$ and $U$ is embedded in $\mathbb{A}^{n+N}$ as $U = V\bigl(r_i(u, v)y_i - q_i(u, v), p_j(u, v)\bigr)$. Since, for each $i$, $x_i$ is the restriction of $y_i$ to $U$, $(x_1, \ldots, x_n)$ is a regular coordinate system if and only if $\det \partial(r_i y_i - q_i, p_j)/\partial(u, v)$ is invertible at every point of $U$; the result follows. □

We will call a Zariski-open subset $U$ of $|M|$ which satisfies the conditions of Lemma 3.3 a *(regular) coordinate chart* with *(regular) coordinates* $x = (x_1, \ldots, x_n)$.

**Definition 3.4.** *Taylor homomorphism.* Let $U$ be a regular coordinate chart in $M$, with coordinates $(x_1, \ldots, x_n)$. For each $a \in U$, there is an injective $\underline{k}$-algebra homomorphism $T_a \colon \mathcal{O}_{M,a} \to \mathbb{F}_a[[X]]$, $X = (X_1, \ldots, X_n)$, which can be described explicitly as follows. Let $p = (p_1, \ldots, p_{N-n})$ (in the notation of (3.1)). By the formal implicit function theorem,

$$p\big(u(a) + X,\, v(a) + V\big) \;=\; U(X, V)\big(V - \varphi(X)\big)\,,$$

where $\varphi(X) \in \mathbb{F}_a[[X]]^{N-n}$, $\varphi(0) = 0$, and $U(X, V)$ is an invertible $(N-n) \times (N-n)$ matrix with entries in $\mathbb{F}_a[[X, V]]$. Let $f \in \mathcal{O}_{M,a}$. Then $f$ is induced by an element $F \in \underline{k}[u, v]_{(a)}$, and $(T_a f)(X) = F\big(u(a) + X,\, v(a) + \varphi(X)\big)$.

The *Taylor homomorphism* $T_a$ induces an isomorphism $\widehat{\mathcal{O}}_{M,a} \to \mathbb{F}_a[[X]]$. Let $D^\alpha \colon \mathbb{F}_a[[X]] \to \mathbb{F}_a[[X]]$ denote the formal derivative $\partial^{|\alpha|}/\partial X^\alpha = \partial^{\alpha_1 + \cdots + \alpha_n}/\partial X_1^{\alpha_1} \cdots \partial X_n^{\alpha_n}$, $\alpha = (\alpha_1, \ldots, \alpha_n) \in \mathbb{N}^n$.

**Lemma 3.5.** *Let $U$ be a regular coordinate chart in $M$, with coordinates $x = (x_1, \ldots, x_n)$. Let $x \in \mathbb{N}^n$. If $f \in \mathcal{O}(U)$, then there is (a unique) $f_\alpha \in \mathcal{O}(U)$ such that, for all $a \in U$,*

$$D^\alpha(T_a f)(X) \;=\; (T_a f_\alpha)(X)\,.$$

*(We will write $f_\alpha = \partial^{|\alpha|} f/\partial x^\alpha$.) More precisely, if $\alpha = (j)$ for some $j$ (where $(j)$ denotes the multiindex with $1$ in the $j$'th place and $0$ elsewhere; i.e., $D^\alpha = \partial/\partial X_j$) and if $f$ is induced by $F = q/r$, where $q(u,v)$, $r(u,v) \in \underline{k}[u,v]$ (in the notation of (3.1)), then $f_{(j)}$ is induced by*

$$F_{(j)} \;=\; \det \frac{\partial(F, p_1, \ldots, p_{N-n})}{\partial(u_j, v_1, \ldots, v_{N-n})} \bigg/ \det \frac{\partial(p_1, \ldots, p_{N-n})}{\partial(v_1, \ldots, v_{N-n})}\,.$$

*Proof.* It suffices to consider the case that $|\alpha| = 1$; i.e., $\alpha = (j)$, for some $j$. Let $a \in U$; say $\big(u(a), v(a)\big) = (0,0)$. ¿From $(T_a f)(X) = F\big(X, \varphi(X)\big)$ (as in Definition 3.4) and from $p\big(X, \varphi(X)\big) = 0$, we obtain

$$\frac{\partial T_a f}{\partial X_j} \;=\; \frac{\partial F}{\partial u_j}\big(X, \varphi(X)\big) + \frac{\partial F}{\partial v}\big(X, \varphi(X)\big) \cdot \frac{\partial \varphi}{\partial X_j}\,,$$

$$0 \;=\; \frac{\partial p}{\partial u_j}\big(X, \varphi(X)\big) + \frac{\partial p}{\partial v}\big(X, \varphi(X)\big) \cdot \frac{\partial \varphi}{\partial X_j}\,.$$

Thus, $\partial T_a f/\partial X_j = F_{(j)}\big(X, \varphi(X)\big)$, where

$$F_{(j)} \;=\; \frac{\left(\det \dfrac{\partial p}{\partial v}\right) \dfrac{\partial F}{\partial u_j} - \dfrac{\partial F}{\partial v} \cdot \left(\dfrac{\partial p}{\partial v}\right)^{\#} \cdot \dfrac{\partial p}{\partial u_j}}{\det\left(\dfrac{\partial p}{\partial v}\right)}\,.$$

($A^{\#}$ means the matrix such that $A \cdot A^{\#} = \det A \cdot I$.) The numerator in this expression equals

$$\det \begin{pmatrix} \dfrac{\partial F}{\partial u_j} & \dfrac{\partial F}{\partial v} \\ \dfrac{\partial p}{\partial u_j} & \dfrac{\partial p}{\partial v} \end{pmatrix}\,,$$

as required.	□

*Remark 3.6.* It is Lemma 3.5 that will show that, when we compute our invariant $\mathrm{inv}_X$ at an $\mathbb{F}$-rational point $a$ of a scheme of finite type over $\underline{k}$, then $\{x : \mathrm{inv}_X(x) \geq \mathrm{inv}_X(a)\}$ and therefore the centre of our blowing-up are nevertheless defined over $\underline{k}$.

*Remark 3.7.* Let $a \in U$. Suppose that $x_i(a) = 0$, $i = 1, \ldots, n$. (We use the notation above.) If $f \in \mathcal{O}(U)$ and $d \in \mathbb{N}$, then the Taylor expansion $(T_a f)(X)$ with respect to the regular coordinate system $x = (x_1, \ldots, x_n)$ can be written in a unique fashion as

$$(T_a f)(X) = c_0(\widetilde{X}) + c_1(\widetilde{X})X_n + \cdots + c_{d-1}(\widetilde{X})X_n^{d-1} + c_d(X)X_n^d \,,$$

where $\widetilde{X} = (X_1, \ldots, X_{n-1})$. Of course, $\widetilde{x} = (x_1, \ldots, x_{n-1})$ forms a regular coordinate system on $N = V(x_n)$ and, for each $q = 0, \ldots, d-1$, $c_q(\widetilde{X})$ is the Taylor expansion at $a$ of the regular function on $N$ given by the restriction of $\dfrac{1}{q!}\dfrac{\partial^q f}{\partial x_n^q}$. Since the Taylor homomorphism is injective, we will write

$$f(x) = c_0(\widetilde{x}) + \cdots + c_{d-1}(\widetilde{x})x_n^{d-1} + c_d(x)x_n^d$$

for the Taylor expansion above, and we will identify each $c_q(\widetilde{x})$, when convenient, with the element of $\mathcal{O}_{N,a}$ induced by $\dfrac{1}{q!}\dfrac{\partial^q f}{\partial x_n^q}$. (In the case of analytic spaces, the preceding expression is just the usual convergent expansion with respect to $x_n$.)

*Properties of the category of spaces*

Let $\mathcal{A}$ denote any of the (algebraic or analytic) categories of local-ringed spaces over $\underline{k}$ listed in (0.2) (1) and (2). Then $\mathcal{A}$ has the following essential features:

(3.8) (1) Let $X \in \mathcal{A}$. If $Y$ is an open or a closed subspace of $X$, then $Y \in \mathcal{A}$. Locally, $X$ is a closed subspace of a manifold $M \in \mathcal{A}$, where:

(2) A *manifold* $M = (|M|, \mathcal{O}_M)$ is a smooth space such that $|M|$ has a neighbourhood basis given by (the supports of) regular coordinate charts as in (0.3). (It follows that if $X$ is a smooth subspace of a manifold $M$, then $X$ is a manifold and is locally a coordinate subspace of a coordinate chart for $M$. In particular, every smooth space $X$ is a manifold and is, therefore, locally pure-dimensional.)

(3) Let $X \in \mathcal{A}$. Then $\mathcal{O}_X$ is a coherent sheaf of rings and $X$ is *locally Noetherian* (i.e., each point of $|X|$ admits an open neighbourhood on which every decreasing sequence of closed subspaces of $X$ stabilizes).

(4) $\mathcal{A}$ is closed under blowing-up. (It follows that if $M \in \mathcal{A}$ is smooth, then a blowing-up $\sigma \colon M' \to M$ with smooth centre $C \subset M$ can be described locally as a quadratic transformation in regular coordinate charts.

We recall that $\mathcal{O}_X$ is a coherent sheaf of rings if and only if every ideal of finite type in $\mathcal{O}_X$ is coherent. "Blowing-up" in (4) can be understood in terms of the universal mapping definition of Grothendieck (cf. [H1, Ch. 0, §2]). We do not need this definition (and therefore do not recall it); it follows from (3) above that if $X$ is a closed subspace of a manifold $M$, then a blowing-up of $X$ is given by the *strict transform* of $X$ by a blowing-up of $M$. (See "Blowing up" and "The strict transform" below.)

*The Zariski topology*

Let $\mathcal{A}$ denote a category of local-ringed spaces over $\underline{k}$ (as in (0.2), for example). Let $X = (|X|, \mathcal{O}_X) \in \mathcal{A}$.

*Definitions and remarks 3.9.* A subset $S$ of $|X|$ will be called a *Zariski-closed subset* of $|X|$ (or of $X$) if $S$ is the support of a closed subspace of $X$ (in $\mathcal{A}$). Suppose that $S$ and $T$ are Zariski-closed subsets of $|X|$; say $S = \operatorname{supp} \mathcal{O}_X/\mathcal{I}$ and $T = \operatorname{supp} \mathcal{O}_X/\mathcal{J}$, where $\mathcal{I}, \mathcal{J} \subset \mathcal{O}_X$ are ideals of finite type that define closed subspaces in $\mathcal{A}$. Then $S \cap T = \operatorname{supp} \mathcal{O}_X/(\mathcal{I}+\mathcal{J})$ and $S \cup T = \operatorname{supp} \mathcal{O}_X/\mathcal{I} \cdot \mathcal{J}$ are Zariski-closed. A *Zariski-open subset* of $|X|$ (or of $X$) means the complement of a Zariski-closed subset. The Zariski-open subsets of $|X|$ are the open sets of the *Zariski topology*. In general, the (original) topology of $|X|$ might be bigger than the Zariski topology. (For example, in the case of analytic spaces.)

We say that $X$ is *Noetherian* if every decreasing sequence of closed subspaces of $X$ (in $\mathcal{A}$) stabilizes. We say that $|X|$ is *Noetherian* if it is Noetherian as a topological space with the Zariski topology; i.e., every decreasing sequence of Zariski-closed subsets stabilizes. If $X$ is Noetherian, then $|X|$ is Noetherian. We say that $X$ (respectively, $|X|$) is *locally Noetherian* if every point of $|X|$ admits an open neighbourhood $U$ (where $U$ is the support of an open subspace in $\mathcal{A}$) such that every decreasing sequence of closed subspaces of $X$ (respectively, Zariski-closed subsets of $|X|$) stabilizes on $U$. Clearly, if $X$ (respectively, $|X|$) is locally Noetherian and $|X|$ is quasi-compact, then $X$ (respectively, $|X|$) is Noetherian. (A real- or complex-analytic space $X$ is Noetherian if and only if $|X|$ is compact.) If $X$ is locally Noetherian, then the intersection of any family of closed subspaces of $X$ is a subspace; hence the intersection of any family of Zariski-closed subsets of $|X|$ is Zariski-closed.

**Lemma 3.10.** *Suppose that $|X|$ is Noetherian. Let $\Sigma$ be a partially ordered set with the property that every decreasing sequence $\sigma_1 \geq \sigma_2 \geq \cdots$ of elements of $\Sigma$ stabilizes. Let $\tau: |X| \to \Sigma$. Then the following are equivalent:*

*(1) $\tau$ is upper-semicontinuous in the Zariski topology; i.e., each $a \in |X|$ admits a Zariski-open neighbourhood $U$ such that $\tau(x) \leq \tau(a)$ for all $x \in U$.*

*(2) $\tau$ takes only finitely many values and, for all $\sigma \in \Sigma$, $S_\sigma := \{x \in |X| : \tau(x) \geq \sigma\}$ is Zariski-closed.*

*Proof.* Assume (1). Let $\sigma \in \Sigma$. Set $W = |X| \backslash S_\sigma$. If $a \in W$, then $a$ admits a Zariski-open neighbourhood $U_a$ in which $\tau(x) \leq \tau(a)$ (and therefore $\tau(x) \not\geq \sigma$); in particular, $U_a \subset W$. Thus $W = \bigcup_{a \in W} U_a$. Since $|X|$ is Noetherian, $W$ is the union of finitely many $U_a$, so that $W$ is Zariski-open, as required. It follows from the hypothesis on $\Sigma$ that $\tau$ takes only finitely many values.

Conversely, assume (2). Let $a \in |X|$. Set $U = \{x \in |X| : \tau(x) \leq \tau(a)\}$. Then $U$ is the complement of the finite union $\bigcup_{\sigma \not\leq \tau(a)} S_\sigma$. □

**Definition 3.11.** *Let $\Sigma$ denote a partially-ordered set. A function $\tau: |X| \to \Sigma$ is Zariski-semicontinuous if: (1) Locally, $\tau$ takes only finitely many values (locally with respect to subspaces of $X$ in $\mathcal{A}$). (2) For all $\sigma \in \Sigma$, $\{x \in |X| : \tau(x) \geq \sigma\}$ is Zariski-closed.*

The Hilbert-Samuel function $H_{X,\cdot}$ and therefore our invariant $\mathrm{inv}_X$ take values in partially-ordered sets satisfying the stabilization hypothesis of Lemma 3.10 (by [BM4, Th. 5.2.1]; cf. Theorem 1.14). Our constructive definition of $\mathrm{inv}_X$ will show that, when the given topology of $|X|$ differs from the Zariski topology, $\mathrm{inv}_X$ is semicontinuous in a sense that is (at least, *a priori*) weaker than that defined above: (In the notation of Theorem 1.14), every point of $|M_j|$ admits a coordinate neighbourhood $U$ with the property that, for all $a \in U$, $V_a := \{x \in U : \mathrm{inv}_X(x) \leq \mathrm{inv}_X(a)\}$ is Zariski-open in $|M_j|U|$. If $|M_j|$ is locally Noetherian (or in the case of a hypersurface $X$), then, as in Lemma 3.10, there is a covering by coordinate charts $U$ in each of which $\mathrm{inv}_X$ takes only finitely many values and (for any value $\iota$), $\{x \in U : \mathrm{inv}_X(x) \geq \iota\}$ is Zariski-closed in $|M_j|U|$.

Of course, if $S \subset |X|$ and $|X|$ is covered by open subsets $U$ such that each $S \cap U$ is the support of a smooth subspace of $X|U$, then $S$ is globally the support of a smooth subspace of $X$. As a consequence, the centres of the blowings-up prescribed by our desingularization algorithm are always smooth *spaces*. Moreover, in the case of analytic spaces (for example), it follows from invariance of $\mathrm{inv}_X$ with respect to finite extension of the base field $\underline{k}$ that $\mathrm{inv}_X$ is actually Zariski-semicontinuous in the stronger sense.

*Blowing-up*

Let $\mathcal{A}$ denote a category of local-ringed spaces over $\underline{k}$ as in (3.8). Let $M = (|M|, \mathcal{O}_M)$ denote a smooth space in $\mathcal{A}$, and let $C$ be a smooth subspace of $M$. Then $C$ is covered by regular coordinate charts $U$ of $M$, each of which admits a coordinate system $x = (w, z)$, $w = (w_1, \ldots, w_{n-r})$, $z = (z_1, \ldots, z_r)$, with respect to which $C \cap U = V(z) = V(z_1, \ldots, z_r)$.

Let $\sigma: M' = \mathrm{Bl}_C M \to M$ denote the blowing-up of $M$ with centre $C$. Then $\sigma$ can be described in local coordinates as follows. Let $U$ denote a regular coordinate chart as above, and let $U' = \sigma^{-1}(U)$. Then

$$U' \cong \{(a, \xi) \in U \times \mathbb{P}^{r-1} : z(a) \in \xi\},$$

where $\mathbb{P}^{r-1}$ denotes the $(r-1)$-dimensional projective space of lines $\xi$ through the origin in $\underline{k}^r$ (or in $\mathbb{A}^r$); if we write $\xi \in \mathbb{P}^{r-1}$ as $\xi = [\xi_1, \ldots, \xi_r]$ in homogeneous coordinates, then

$$U' = \{(a, \xi) \in U \times \mathbb{P}^{r-1} : z_i(a)\xi_j = z_j(a)\xi_i,\ 1 \leq i,\ j \leq r\}.$$

Therefore, $U' = \bigcup_{i=1}^{r} U'_i$, where, for each $i$,

$$U'_i = \{(a, \xi) \in U' : \xi_i = 1\}.$$

It follows that, for each $i$, $U'_i$ is a regular coordinate chart with coordinates $x' = (w', z')$, $w' = (w'_1, \ldots, w'_{n-r})$, $z' = (z'_1, \ldots, z'_r)$ given by

$$w'(a, \xi) = w(a)$$
$$z'_j(a, \xi) = \begin{cases} z_i(a), & j = i, \\ \xi_j = z_j(a)/z_i(a), & j \neq i. \end{cases}$$

In particular, suppose that $f \in \mathcal{O}_{M,a}$, where $a \in U$ and $w(a) = 0$, $z(a) = 0$; if $a' \in \sigma^{-1}(a) \cap U'_i$, then the Taylor expansion of $f \circ \sigma$ at $a'$ is given by formal substitution of $w = w'$, $z_i = z'_i$ and $z_j = z'_i\big(z'_j(a') + z'_j\big)$, $j \neq i$, in the Taylor expansion of $f$ at $a$.

*Example 3.12.* Let $M = (|M|, \mathcal{O}_M)$ be a smooth scheme of finite type over $\underline{k}$, and let $U$ be a regular coordinate chart with coordinates $x = (x_1, \ldots, x_n)$ as in (3.1). As before, suppose that $x = (w, z)$ such that $C \cap U = V(z)$. Consider $U_i'$, say for $i = 1$. (Using the notation of (3.1)) we have a commutative diagram

$$\begin{array}{ccc} U_1' & \hookrightarrow & \mathbb{A}^{N+(r-1)} \\ & \searrow & \downarrow \\ & & \mathbb{A}^n \end{array}$$

where $U_1'$ is embedded in $\mathbb{A}^{N+(r-1)}$ as $V(p_1, \ldots, p_{N-n}, q_j := u_{n-r+j} - u_{n-r+1}\xi_j, j = 2, \ldots, r)$ (with respect to the affine coordinates $(u, v, \xi_2, \ldots, \xi_r)$ of $\mathbb{A}^{N+(r-1)}$) and the vertical projection is given by $(u_1, \ldots, u_{n-r}, u_{n-r+1}, \xi_2, \ldots, \xi_r)$. $U_1' \to \underline{k}^n$ is an étale covering since

$$\det \frac{\partial(p, q_2, \ldots, q_r)}{\partial(v, u_{n-r+2}, \ldots, u_n)} \ = \ \det \frac{\partial p}{\partial v} \ .$$

*The strict transform*

We continue to use the notation of the preceding subsection. Let $\sigma \colon M' \to M$ denote a blowing-up with smooth centre $C \subset M$, and let $H = \sigma^{-1}(C)$. Let $X$ be a closed subspace of $M$. First suppose that $X$ is a *hypersurface*; i.e., $\mathcal{I}_X$ is principal. Let $a \in M$ and let $f \in \mathcal{I}_{X,a}$ denote a generator of $\mathcal{I}_{X,a}$. If $a' \in \sigma^{-1}(a)$, then we define $\mathcal{I}_{X',a'}$ as the principal ideal in $\mathcal{O}_{M',a'}$ generated by $f' = y_{\text{exc}}^{-d} f \circ \sigma$, where $y_{\text{exc}}$ denotes a generator of $\mathcal{I}_{H,a'}$ and $d = \mu_{C,a}(f)$. (Thus $d$ is the largest power of $y_{\text{exc}}$ to which $f \circ \sigma$ is divisible in $\mathcal{O}_{M',a'}$.) In this way we obtain a coherent sheaf of principal ideals $\mathcal{I}_{X'}$ in $\mathcal{O}_{M'}$; the *strict transform $X'$ of $X$ by $\sigma$* means the corresponding closed subspace of $M'$.

In local coordinates as above, suppose that $w(a) = 0$, $z(a) = 0$, and let $a' \in U_1'$. Then $y_{\text{exc}} = z_1'$ and (the Taylor expansion at $a'$ of) $f'$ is given by

$$f'(w', z') \ = \ (z_1')^{-d} f\bigl(w', z_1', z_1'(\widetilde{z}'(a) + \widetilde{z}')\bigr) \, ,$$

where $\widetilde{z}' = (z_2', \ldots, z_r')$. We will also call $f'$ the "strict transform" of $f$ by $\sigma$, although $f'$ is, of course, only defined up to multiplication by an invertible factor.

Now consider an arbitrary closed subspace $X$ of $M$. Then the *strict transform of $X$ by $\sigma$* can be defined as the closed subspace $X'$ of $M'$ such that, locally at each $a' \in M'$, $X'$ is the intersection of the strict transforms of all hypersurfaces containing $X$ near $a = \sigma(a')$. To be precise: If $a' \in M'$, let $\mathcal{I}_{X',a'} \subset \mathcal{O}_{M',a'}$ denote the ideal generated by the strict transforms $f'$ of all $f \in \mathcal{I}_{X,a}$.

**Proposition 3.13.** *Let $a' \in M'$. Then $\mathcal{I}_{X',a'}$ coincides with the ideal $\{f \in \mathcal{O}_{M',a'} : y_{\text{exc}}^k f \in \mathcal{I}_{\sigma^{-1}(X),a'},$ for some $k \in \mathbb{N}\}$. ($\mathcal{I}_{\sigma^{-1}(X),a'}$ is the ideal generated by $\sigma_{a'}^*(\mathcal{I}_{X,\sigma(a')})$.)*

We will give a simple proof of Proposition 3.13 below, as an application of the diagram of initial exponents. By Proposition 3.13, $\mathcal{I}_{X'} = \sum_k [\mathcal{I}_{\sigma^{-1}(X)} : y_{\text{exc}}^k]$, so that $\mathcal{I}_{X'}$ is an ideal of finite type (since $X$ is locally Noetherian).

*Remark 3.14.* (In a category satisfying the conditions (3.8)), we define the *reduced space* $X_{\text{red}}$ corresponding to $X$ using the coherent sheaf of ideals $\mathcal{I}_{X_{\text{red}}} = \sqrt{\mathcal{I}_X}$. It is easy to see that if $X'$ is the strict transform of $X$ by a blowing-up, as above, then $(X')_{\text{red}} = (X_{\text{red}})'$.

*Remark 3.15.* Let $X''$ denote the smallest closed subspace of $\sigma^{-1}(X)$ containing $\sigma^{-1}(X)\backslash H$, where $H = \sigma^{-1}(C)$. ($X''$ exists by local Noetherianness.) Of course, $X'' \subset X'$. In the case of schemes or analytic spaces over an algebraically closed field, $X'' = X'$. (A consequence of Hilbert's Nullstellensatz.) The following trivial example shows that $X'' \neq X'$ in general.

*Example 3.16.* Let $X \subset M = \mathbb{R}^2$ denote the real analytic subspace determined by the coherent sheaf of ideals generated by $x^4(x-1)^2 + y^2$. Let $\sigma: M' \to M$ denote the blowing-up with centre $\{0\}$. Then the strict transform $X' \subset U'$, where $U' \subset M'$ is a coordinate chart in which $\sigma$ is given by $x = u$, $y = uv$. In $U'$, $X'$ is defined by $u^2(u-1)^2 + v^2 = 0$, so that $|X'| = \{(0,0), (1,0)\}$, but $|X''| = \{(1,0)\}$.

*The diagram of initial exponents*

The material in this subsection is needed only in Chapters III and IV (but we will also use the diagram to give a simple proof of Proposition 3.13 and to extend Remark 1.11 to the general case).

Let $\mathbb{K}$ be a field and let $\mathbb{K}[[X]]$ denote the ring of formal power series in $X = (X_1, \ldots, X_n)$. If $\alpha = (\alpha_1, \ldots, \alpha_n) \in \mathbb{N}^n$, put $|\alpha| = \alpha_1 + \cdots + \alpha_n$. The lexicographic ordering of $(n+1)$-tuples $(|\alpha|, \alpha_1, \ldots, \alpha_n)$ induces a total ordering of $\mathbb{N}^n$. Let $F \in \mathbb{K}[[X]]$. Write $F = \sum_{\alpha \in \mathbb{N}^n} F_\alpha X^\alpha$, where $X^\alpha = X_1^{\alpha_1} \cdots X_n^{\alpha_n}$. Let $\operatorname{supp} F = \{\alpha \in \mathbb{N}^n : F_\alpha \neq 0\}$. The *initial exponent* $\exp F$ is defined as the smallest element of $\operatorname{supp} F$. If $\alpha = \exp F$, then $F_\alpha x^\alpha$ is called the initial monomial $\operatorname{mon} F$ of $F$.

The following formal division theorem of Hironaka [H1] is a simple generalization of the Euclidean division algorithm. Let $G^1, \ldots, G^s \in \mathbb{K}[[X]]$, and let $\alpha^i = \exp G^i$, $i = 1, \ldots, s$. We associate to $\alpha^1, \ldots, \alpha^s$ the following decomposition of $\mathbb{N}^n$: Set $\Delta_i = (\alpha^i + \mathbb{N}^n) - \bigcup_{j=1}^{i-1} \Delta_j$, $i = 1, \ldots, s$, and put $\square_0 = \mathbb{N}^n - \bigcup_{i=1}^{s} \Delta_i$. We also define $\square_i \subset \mathbb{N}^n$ by $\Delta_i = \alpha^i + \square_i$, $i = 1, \ldots, s$.

**Theorem 3.17.** *For every $F \in \mathbb{K}[[X]]$, there are unique $Q_i \in \mathbb{K}[[X]]$, $i = 1, \ldots, s$, and $R \in \mathbb{K}[[X]]$ such that $\operatorname{supp} Q_i \subset \square_i$, $i = 1, \ldots, s$, $\operatorname{supp} R \subset \square_0$, and*

$$F = \sum_{i=1}^{s} Q_i G^i + R.$$

*Remark 3.18.* Let $\underline{m}$ denote the maximal ideal of $\mathbb{K}[[X]]$. In Theorem 3.17, if $k \in \mathbb{N}$ and $F \in \underline{m}^k$, then $R \in \underline{m}^k$ and each $Q_i \in \underline{m}^{k-|\alpha^i|}$ (where $\underline{m}^\ell$ means $\mathbb{K}[[X]]$ if $\ell \leq 0$).

Let $I$ be an ideal in $\mathbb{K}[[X]]$. The *diagram of initial exponents* $\mathfrak{N}(I) \subset \mathbb{N}^n$ is defined as

$$\mathfrak{N}(I) = \{\exp F : F \in I\}.$$

Clearly $\mathfrak{N}(I) + \mathbb{N}^n = \mathfrak{N}(I)$. Let $\mathcal{D}(n) = \{\mathfrak{N} \subset \mathbb{N}^n : \mathfrak{N} + \mathbb{N}^n = \mathfrak{N}\}$. If $\mathfrak{N} \in \mathcal{D}(n)$, then there is a smallest finite subset $\mathfrak{V}$ of $\mathfrak{N}$ such that $\mathfrak{N} = \mathfrak{V} + \mathbb{N}^n$; $\mathfrak{V} = \{\alpha \in \mathfrak{N} : \mathfrak{N}\backslash\{\alpha\} \in \mathcal{D}(n)\}$. We call $\mathfrak{V}$ the *vertices* of $\mathfrak{N}$. The following is a simple consequence of Theorem 3.17.

**Corollary 3.19.** *Let $\alpha^i$, $i = 1, \ldots, s$, denote the vertices of $\mathfrak{N}(I)$. Choose $G^i \in I$ such that $\alpha^i = \exp G^i$, $i = 1, \ldots, s$ (we say that $G^i$ represents $\alpha^i$), and let $\{\Delta_i, \Box_0\}$ denote the decomposition of $\mathbb{N}^n$ determined by the $\alpha^i$, as above. Then:*

*(1) $\mathfrak{N}(I) = \bigcup \Delta_i$ and the $G^i$ generate $I$.*

*(2) There is a unique set of generators $F^i$ of $I$, $i = 1, \ldots, s$, such that, for each $i$, $\mathrm{supp}(F^i - x^{\alpha^i}) \subset \Box_0$; in particular, $\mathrm{mon}\, F^i = x^{\alpha^i}$.*

We call $F^1, \ldots, F^s$ the *standard basis* of $I$ (with respect to the given total ordering of $\mathbb{N}^n$). If $\mathfrak{N} \in \mathcal{D}(n)$, let $\mathbb{K}[[X]]^{\mathfrak{N}} = \{F \in \mathbb{K}[[X]] : \mathrm{supp}\, F \cap \mathfrak{N} = \emptyset$; i.e., $\mathrm{supp}\, F \subset \Box_0\}$. Clearly, $\mathbb{K}[[X]]^{\mathfrak{N}}$ *is stable with respect to formal differentiation.*

Now let $H_I$ denote the Hilbert-Samuel function of $\mathbb{K}[[X]]/I$; i.e., $H_I(k) = \dim_{\mathbb{K}} \mathbb{K}[[X]]/(I + \underline{m}^{k+1})$, $k \in \mathbb{N}$. By Remark 3.18 and Corollary 3.19, we have:

**Corollary 3.20.** *For every $k \in \mathbb{N}$, $H_I(k) = \#\{\alpha \in \mathbb{N}^n : \alpha \notin \mathfrak{N}(I)$ and $|\alpha| \leq k\}$.*

It follows from Corollary 3.20 that $H_I(k)$ coincides with a polynomial in $k$, for $k$ sufficiently large (the "Hilbert-Samuel polynomial").

*Remark 3.21.* The preceding definitions make sense and the results above (*except for* Remark 3.18 and Corollary 3.20) hold, more generally, for any total ordering of $\mathbb{N}^n$ which is compatible with addition in the sense that: For any $\alpha, \beta, \gamma \in \mathbb{N}^n$, $\gamma \geq 0$, and $\alpha \leq \beta \Rightarrow \alpha + \gamma \leq \beta + \gamma$.

In order to prove Proposition 3.13, we will use the total ordering of $\mathbb{N}^n$ given by the lexicographic ordering of $(\alpha_1, |\alpha|, \alpha_2, \ldots, \alpha_n)$, $\alpha \in \mathbb{N}^n$. We then have:

**Lemma 3.22.** *Let $I$ be an ideal in $\mathbb{K}[[Y]] = \mathbb{K}[[Y_1, \ldots, Y_n]]$, and let $J$ denote the ideal $J = \{G(Y) \in \mathbb{K}[[Y]] : Y_1^k G(Y) \in I$, for some $k \in \mathbb{N}\}$. Suppose that $F_i(Y) \in I$, $i = 1, \ldots, s$, represent the vertices of $\mathfrak{N}(I)$; for each $i$, write $F_i(Y) = Y_1^{k_i} G_i(Y)$, where $G_i(Y)$ is not divisible by $Y_1$. Then $J$ is generated by the $G_i$.*

*Proof.* This is an immediate consequence of the following variant of Remark 3.18 which holds for the given ordering of $\mathbb{N}^n$: In the formal division algorithm 3.17, if $F \in (Y_1)^k$, then $R \in (Y_1)^k$ and each $Q_i \in (Y_1)^{k-k_i}$. ($(Y_1)$ denotes the ideal generated by $Y_1$.) □

*Proof of Proposition 3.13.* We can choose coordinates at $a = \sigma(a')$ and $a'$ so that $\widehat{\mathcal{O}}_{M,a} \cong \underline{k}[[X_1, \ldots, X_n]]$, $\widehat{\mathcal{O}}_{M',a'} \cong \underline{k}[[Y_1, \ldots, Y_n]]$ and $\widehat{\sigma}^*_{a'} : \widehat{\mathcal{O}}_{M,a} \to \widehat{\mathcal{O}}_{M',a'}$ has the form $X_\ell = Y_\ell$, $\ell = 1, \ldots, q$ (where $q \geq 1$), $X_\ell = Y_1(\eta_\ell + Y_\ell)$, $\ell = q+1, \ldots, n$. Put $I = \widehat{\mathcal{I}}_{\sigma^{-1}(X), a'} \subset \underline{k}[[Y]]$ and $J = \{G(Y) \in \underline{k}[[Y]] : Y_1^k G(Y) \in I$, for some $k \in \mathbb{N}\}$. Suppose that $H_j(X)$, $j = 1, \ldots, r$, generate $\widehat{\mathcal{I}}_{X,a}$. We can find polynomials $P_{ij}(Y) \in \underline{k}[Y]$, $i = 1, \ldots, s$, $j = 1, \ldots, r$, such that the

$$F_i(Y) := \sum_j P_{ij}(Y)(H_j \circ \sigma)(Y) \in I$$

represent the vertices of $\mathfrak{N}(I)$. Each $F_i(Y)$ is the pullback by $\sigma$ of

$$\sum_j P_{ij}\left(X_1, \ldots, X_q, \frac{X_{q+1}}{X_1} - \eta_{q+1}, \ldots, \frac{X_n}{X_1} - \eta_n\right) H_j(X)$$
$$= (X_1)^{-q_i} \sum_j Q_{ij}(X) H_j(X),$$

for some $q_i \in \mathbb{N}$, where the $Q_{ij} \in \underline{k}[X]$. Write $G_i(X) = \sum_j Q_{ij}(X) H_j(X)$, for each $i$. Thus each $F_i(Y) = (Y_1)^{-q_i}(G_i \circ \sigma)(Y)$. We write $(G_i \circ \sigma)(Y) = Y_1^{m_i} G'_i(Y)$, where $G'_i$ is not divisible by $Y_1$, so that $m_i \geq q_i$. Then we have $F_i(Y) = Y_1^{m_i - q_i} G'_i(Y)$, for each $i$, where each $G'_i \in \widehat{\mathcal{I}}_{X',a'}$. But the $G'_i$ generate $J$, by Lemma 3.22. (The preceding formal argument suffices to prove the proposition because, for any ideal $\mathcal{J}$ in $\mathcal{O}_{M',a'}$, $\widehat{\mathcal{J}} \cap \mathcal{O}_{M',a'} = \mathcal{J}$.) $\square$

*Remark 3.23.* The diagram of initial exponents can be used to generalize the geometric definition of $\mathrm{inv}_X$ in year zero, given in Remark 1.11 in the case of a hypersurface. Let $a \in M$ and let $(x_1, \ldots, x_n)$ denote a coordinate system at $a$, so that $\widehat{\mathcal{O}}_{M,a} \cong \underline{k}[[x_1, \ldots, x_n]]$ via the Taylor homomorphism. Let $w = (w_1, \ldots, w_n)$ be an $n$-tuple of positive real numbers ("weights" for the coordinates). For $f(x) = \sum f_\alpha x^\alpha \in \underline{k}[[x]]$, we define the *weighted order* $\mu_w(f) := \min\{\langle w, \alpha \rangle : f_\alpha \neq 0\}$ (where $\langle w, \alpha \rangle := \sum w_i \alpha_i$) and the *weighted initial exponent* $\exp_w(f) := \min\{\alpha : f_\alpha \neq 0\}$, where the multiindexes $\alpha \in \mathbb{N}^n$ are totally ordered using lexicographic ordering of the sequences $(\langle w, \alpha \rangle, \alpha_1, \ldots, \alpha_n)$. Set $I = \widehat{\mathcal{I}}_{X,a}$. Write $\mathfrak{N}_w(I)$ for the (*weighted*) diagram of initial exponents $\{\exp_w(f) : f \in I\}$. If $\mathfrak{N} \in \mathcal{D}(n)$, we define the *essential variables* of $\mathfrak{N}$ as the indeterminates $x_j$ which occur (to positive power) in some monomial $x^\alpha$, where $\alpha \in \mathfrak{V}$ (the vertices of $\mathfrak{N}$).

For the given coordinate system $x = (x_1, \ldots, x_n)$, let $d(x)$ denote the supremum of $n$-tuples $(d_1, \ldots, d_n) \in (\mathbb{Q} \cup \{\infty\})^n$, ordered lexicographically, such that:

(1) $1 = d_1 \leq d_2 \leq \cdots \leq d_n$;

(2) $x_1, \ldots, x_r$ are the essential variables of $\mathfrak{N}(I)$, for some $r$, and $d_1 = \cdots = d_r = 1$;

(3) $\mathfrak{N}(I) = \mathfrak{N}_w(I)$, where $w_j = 1/d_j$, $j = 1, \ldots, n$.

Here $\mathfrak{N}(I)$ denotes the diagram with respect to the standard ordering $(|\alpha|, \alpha_1, \ldots, \alpha_n)$ of $\mathbb{N}^n$. Conditions (2) and (3) above are together equivalent to: If $F^1, \ldots, F^s$ is the standard basis of $I$ (with respect to the standard ordering of $\mathbb{N}^n$), then $\mu(F^i) = \mu_w(F^i)$ for each $i$. ($\mu(F^i)$ is the standard order $\mu_a(F^i)$.)

Set $d = \sup d(x)$ (sup over all coordinate systems), $d = (d_1, \ldots, d_n)$. Then

$$\mathrm{inv}_X(a) = \left( H_{X,a}, 0; \frac{d_2}{d_1}, 0; \ldots; \frac{d_t}{d_{t-1}}, 0; \infty \right),$$

where $d_t$ is the last finite $d_i$. Moreover, beginning with any coordinate system, we can make an explicit change of variables to obtain coordinates $x = (x_1, \ldots, x_n)$ satisfying a criterion which guarantees that $d(x) = d$; in these coordinates, the centre of the first blowing-up in the resolution algorithm is $x_1 = \cdots = x_t = 0$. Again, there is a correspondence between the weighted initial ideals of $I$ with respect to two coordinate systems that realize $d$, analogous to that of Remark 1.11.